\documentclass{article}
\usepackage{amsmath}
\usepackage{amssymb}
\usepackage{graphicx}
\usepackage[urlcolor=blue,citecolor=red,colorlinks=true]{hyperref}
\usepackage{paralist}
\usepackage{subcaption}

\newcommand{\ep}{\varepsilon}

\title{Manifold Learning for Accelerating\\Coarse-Grained Optimization}
\author{Dmitry Pozharskiy, Noah J. Wichrowski, Andrew B. Duncan,\\ Grigorios A. Pavliotis, and Ioannis G. Kevrekidis}
\date{}

\begin{document}
	\maketitle
	
	\begin{abstract}
		\noindent Algorithms proposed for solving high-dimensional optimization problems with no derivative information frequently encounter the ``curse of dimensionality,''  becoming ineffective as the dimension of the parameter space grows. One feature of a subclass of such problems that are {\em effectively} low-dimensional is that only a few parameters (or combinations thereof) are important for the optimization and must be explored in detail. Knowing these parameters/combinations \emph{in advance} would greatly simplify the problem and its solution. We propose the data-driven construction of an \emph{effective} (coarse-grained, ``trend'') optimizer, based on data obtained from ensembles of brief simulation bursts with an ``inner" optimization algorithm, that has the potential to accelerate the exploration of the parameter space. The trajectories of this ``effective optimizer'' quickly become attracted onto a slow manifold parameterized by the few relevant parameter combinations. We obtain the parameterization of this low-dimensional, effective optimization manifold \emph{on the fly} using data mining/manifold learning techniques on the results of simulation (inner optimizer iteration) burst ensembles and exploit it locally to ``jump'' forward along this manifold. As a result, we can bias the exploration of the parameter space towards the few, important directions and, through this ``wrapper algorithm,'' speed up the convergence of traditional optimization algorithms.
	\end{abstract}
	
	\section{Introduction}
	The design of complex engineering systems often leads to high-dimensional optimization problems with computationally expensive objective function evaluations, often given in the form of a (computational) black-box. In such cases the derivative information may be unavailable or impractical to obtain in closed form, too expensive to compute numerically, or unreliable to estimate if the objective function is noisy. These difficulties may render derivative-based optimization methods impractical for such problems, and so-called derivative-free methods must be used.
	
	The first such derivative-free algorithms appeared quite early: the direct search method~\cite{hooke1961direct} and the Nelder-Mead algorithm~\cite{nelder1965simplex}. Since then a variety of algorithms have been proposed, including trust-region methods~\cite{conn2000trust}; deterministic global algorithms~\cite{jones1993lipschitzian,huyer1999global}; algorithms utilizing surrogate models~\cite{barton1994metamodeling,jones2001taxonomy}; and stochastic global methods such as genetic algorithms~\cite{holland1992adaptation}, simulated annealing~\cite{kirkpatrick1983optimization} and particle swarm optimization~\cite{eberhart1995new}. A broader overview of the aforementioned methods along with additional ones used for high-dimensional problems can be found at~\cite{conn2009introduction,shan2010survey,rios2013derivative}.
	
	Many of these methods have been applied successfully to low-dimensional problems where derivative information is not available; once the dimension of the parameter space grows, however, they run into the ``curse of dimensionality,'' where the required sampling of the parameter space grows exponentially or the convergence becomes too slow. There is case-dependent evidence that, for certain classes of problems, out of the vast parameter space only a few parameters or combinations of parameters suffice to describe most of the variance in the objective function values, with the rest having little effect~\cite{li2001high,montgomery2006general,montgomery2007general,roslund2014dynamic}. It is observations of this nature that we aim to exploit in our proposed method, borrowing additional ideas from the field of fast/slow (singularly perturbed) dynamical systems and data mining/manifold learning techniques.
	
	It has been observed that complex {\em dynamical systems} such as molecular dynamics (MD) simulations or complex reaction network dynamics may possess a low-dimensional, attracting, ``slow'' manifold. The dynamics of the system after being initialized at a random state quickly approach the slow manifold and then evolve ``along it'' (close to it). A reduced model of the system in terms of the slow variables parameterizing this manifold would greatly simplify the understanding of the system's behavior (and its computation). However, such a model is often unavailable in closed form. In previous work~\cite{kevrekidis2003equation} we have shown how short ``bursts'' of a microscopic simulator can evolve the system close to and then ``along'' an
	underlying slow manifold. Essentially, after the short burst is attracted to the slow manifold, we can \emph{observe} the evolution on a \emph{restricted} set of coarse-grained observables that parameterize the slow manifold when these coarse variables are known \emph{a priori}. We can then perform a ``large'' time step by extrapolating the few macroscopic variable values and \emph{lifting} the new state back into the full space to initialize a new set of computation bursts for the microscopic simulator. This can achieve significant acceleration of the \emph{effective} complex system dynamics. If the macroscopic variables are not known, then a reduced description of the manifold can be derived \emph{on the fly} by using data-driven dimensionality reduction techniques, which uncover the few intrinsic variables that are adequate to describe a high-dimensional data set locally.
	
	The high-dimensional optimization problem can be treated in the same vein by making two assumptions: (a) we have an ``inner optimizer'', analogous to a microscopic simulator, that samples the parameter space and produces a series of ``pseudo''-Langevin trajectories, (b) we postulate that there exists an attracting, slow manifold which can be parameterized in terms of a few parameters or combinations thereof, and the inner optimizer is quickly attracted to it. In the following sections we will show that the trajectory produced by a specific version of simulated annealing (SA){\textemdash}our ``inner optimizer''{\textemdash}at constant temperature can be described by an effective Stochastic Differential Equation (SDE) whose drift contains information about the local gradient of the objective function. To be precise, since we do not vary the temperature as the optimization progresses, our "inner optimizer" is a random walk Metropolis-Hastings (RWMH) algorithm.  Running short bursts of this constant temperature SA, we create ``pseudo'' dynamics that can be thought of as the (approximate) dynamics of an actual dynamical system. After initializing at a random point in parameter space, the algorithm is quickly attracted to the low-dimensional manifold, and by applying either linear or nonlinear dimensionality reduction techniques we can obtain a useful local parameterization of this manifold. We can estimate the drift of the effective SDE using established parametric inference methods and thus estimate the local \emph{effective} gradient of the objective function \cite{schillings2017analysis}. This can be used subsequently in an algorithm such as gradient descent in a reduced parameter space. The new point is lifted back to full space, using local Principal Component Analysis or geometric harmonics~\cite{coifman2006geometric}, and the entire procedure is repeated, leading to an acceleration of the overall optimization.
	
	\section{Methods}
	\subsection{Inner optimization loop}
	The Langevin equation was introduced as a stochastic global optimization method shortly after the first appearance of simulated annealing~\cite{gidas1985global,geman1986diffusions}. It is a gradient descent method with the addition of a ``thermal'' noise that allows the trajectories to escape local minima and thus enhance their ability to explore the parameter space. However, it may become impractical for the problems we are considering since, as we discussed above, the gradient information is explicitly unavailable. The equation reads
	\begin{equation} \label{eq:langevin}
		dx_t=-\nabla f(x_t)\,dt+\sqrt{2T}\,dW_t,
	\end{equation}
	where $x_t\in\mathbb{R}^n$, $f$ is the objective function, $W_t$ is an $n$-dimensional standard Brownian motion, and $T$ is the temperature parameter. It can be shown that under an appropriate temperature schedule $T(t)$, the algorithm converges
	weakly to the global minimum \cite{geman1986diffusions}. The equilibrium distribution is the Gibbs distribution, with density
	\begin{equation*}
		\pi(x;\,T)=\dfrac{\exp\left[-\frac{f(x)}{T}\right]}{\int_{\mathbb{R}^n}\exp\left[-\frac{f(y)}{T}\right]\,dy}.
	\end{equation*}
	The simulated annealing algorithm admits the same equilibrium distribution at an equal $T$ and can be viewed as an adaptation of the Metropolis-Hastings algorithm~\cite{metropolis1953equation} with time dependent acceptance probability due to the temperature schedule. The acceptance probability is given by
	\begin{equation*}
		a=\min\left(1,\,\exp\left[\frac{-\big(f(y)-f(x)\big)}{T}\right]\right),
	\end{equation*}
	where $x$ is the current point and $y$ is the new trial point that comes from a symmetric proposal density $g(y|x)$. Hence, better points are always accepted and worse points are accepted with probability $0<a<1$, which is greater at higher temperatures $T$.
	
\begin{figure}[htb]
		\centering
		\begin{subfigure}[]{0.95\linewidth}
			\includegraphics[width=\linewidth]{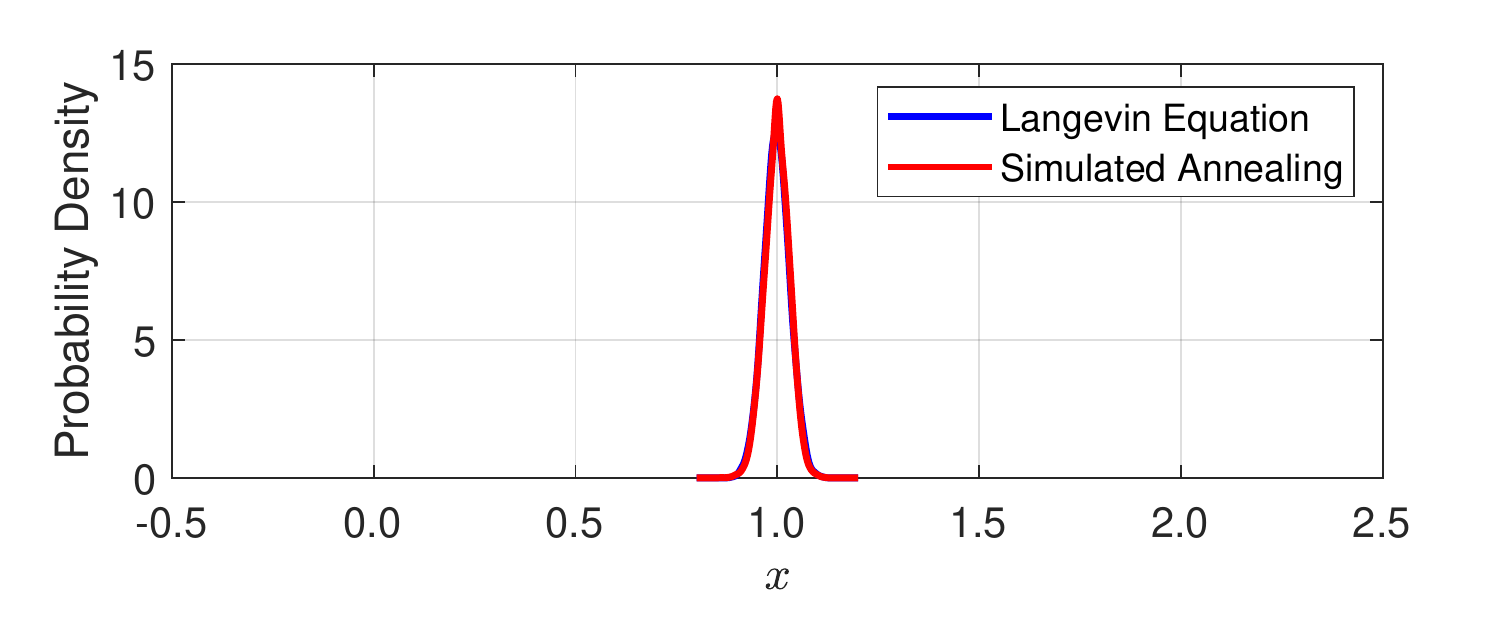}
			\caption{After one step ($t=0.001$).}
			\label{fig:onestep}
		\end{subfigure}
		\begin{subfigure}[]{0.95\linewidth}
			\includegraphics[width=\linewidth]{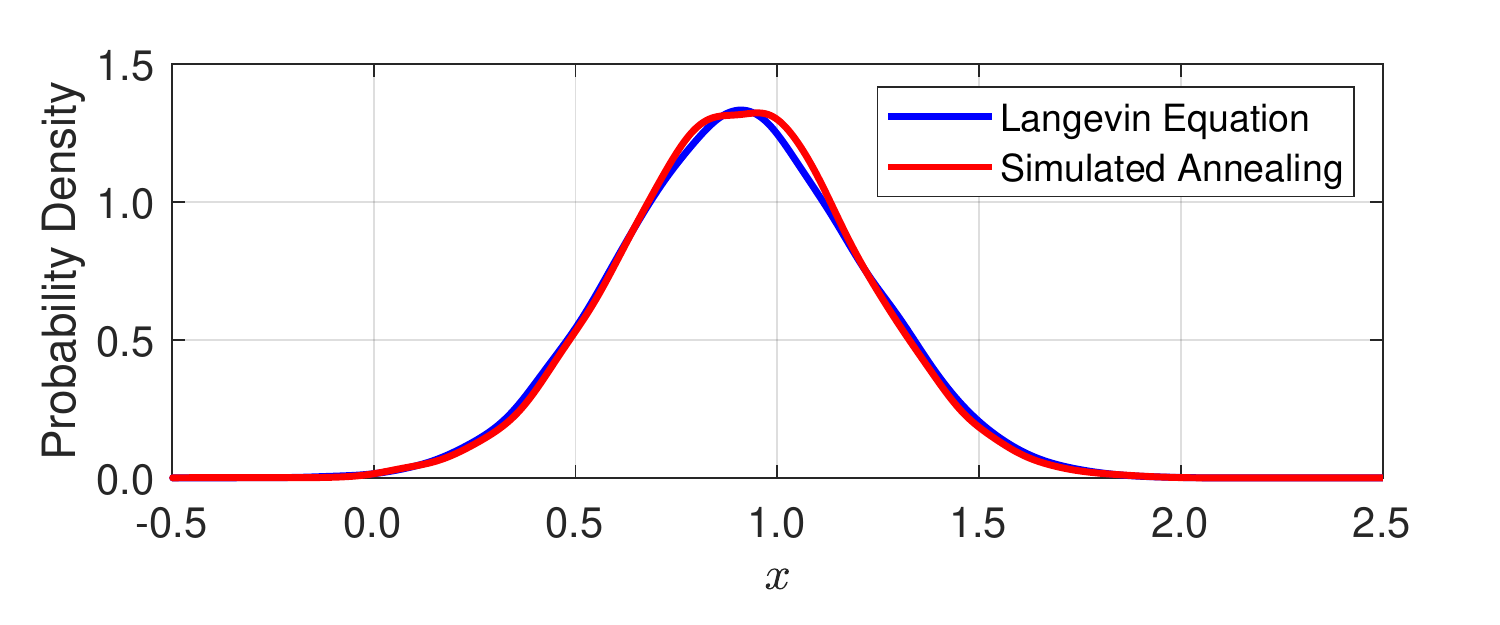}
			\caption{After 100 steps ($t=0.1$).}
			\label{fig:100step}
		\end{subfigure}
		\caption[Probability density of points using Simulated Annealing and the Langevin equation]{Evolution in time of the probability density of current points using either the Langevin equation or SA at constant $T$. The objective function is $f(x)=0.5x^2$; $10^4$ realizations are used with starting point $x=1$, $T=0.5$; and the time step is $dt=10^{-3}$.}
		\label{fig:steps}
	\end{figure}

	Consider now the simple, one-dimensional case at constant temperature $T$ (a RWMH protocol) using proposal density $g(y|x)=\mathcal{N}(x,\,2T\,\delta t)$ and the acceptance probability defined above. It can be shown \cite{duncan2017nonreversible,fathi2017improving} that, at the limit of \emph{small time steps} $\delta t$ or \emph{large temperatures} $T$, the density of accepted points after one step converges to a normal distribution $\mathcal{N}(x-f'(x)\,\delta t,\,2T\,\delta t)$, which corresponds to the density of a new point using an Euler-Maruyama discretization of the Langevin equation \cite{locatelli2002simulated}. Figure~\ref{fig:steps} shows the density of sample points after one step and after 100 steps using both algorithms. The two distributions visually almost coincide.
	
	Using the above procedure we can obtain ``pseudo'' trajectories in the parameter space that are analogous to the trajectories produced by the Langevin equation and contain information about an {\em effective} gradient of the objective function without explicitly computing it.
	
	\subsection{Dimensionality reduction}
	In our previous discussion on dynamical systems, we mentioned that the long-term dynamics of the full system can often be usefully restricted to the dynamics of a few {\em slow} variables. These variables can be a collection of macroscopic variables that are available from our intimate knowledge of the system, or they can be estimated ``on the fly'' using dimension reduction techniques. Such techniques can be applied to large, high-dimensional data sets to uncover the few {\em intrinsic} variables that are sufficient to describe most of the system's long-term behavior. We will use the latter approach in our method, since it can be quite challenging to identify beforehand the few parameters that are important to the optimization.
	
	One of the most common methods for dimension reduction is Principal Component Analysis (PCA)~\cite{jolliffe1986principal}. It tries to identify a hyperplane that best fits the data by finding an orthogonal basis, where the first vector points in the direction of maximum variance in the data set and all subsequent vectors maximize variance in orthogonal directions. The basis vectors are called \emph{Principal Components}, and they can be found by an eigenvalue decomposition of the covariance matrix of the data set after it has been centered. If the eigenvalues are sorted and the relationship $\lambda_1>\lambda_2>\cdots>\lambda_k\gg\lambda_{k+1}>\cdots>\lambda_n$ holds ($n$ is the dimension of the original space), then there is a gap in the eigenvalue spectrum and we can reduce the dimensionality of our data set by projecting it onto the first $k$ principal components. PCA is a well-documented technique, but its major limitation is that it can parameterize only linear manifolds.
	
	Nonlinear manifold learning techniques are required if the data lie on a curved manifold. One such method is Diffusion Maps (DMaps)~\cite{coifman2005geometric}. For a data set of size $m$, the algorithm starts by constructing a weight matrix $W\in\mathbb{R}^{m\times m}$:
	\begin{equation*}
		W_{ij}=\exp\left(\frac{-\lVert x_i-x_j\rVert^2}{\ep^2}\right),\;i,j=1,\ldots,m,
	\end{equation*}
	where $x_i,x_j\in\mathbb{R}^n$, $\lVert\cdot\rVert$ is an appropriate norm, and $\ep$ is a characteristic distance between data points. Next, we construct the diagonal matrix $D\in\mathbb{R}^{m\times m}$ with $D_{ii}=\sum_j W_{ij}$ and compute $\tilde{W}=D^{-\alpha}WD^{-\alpha}$, where $0\leq\alpha\leq1$ is a normalization parameter. Then, we construct the diagonal matrix $\tilde{D}\in\mathbb{R}^{m\times m}$ with $\tilde{D}_{ii}=\sum_j\tilde{W}_{ij}$ and compute the row-stochastic matrix $K=\tilde{D}^{-1}\tilde{W}$. The matrix $K$ is the transition probability matrix of a Markov chain defined on the data set, whose states are the individual data points.
	
	The eigenvectors $\psi_0,\psi_1,\ldots,\psi_{m-1}$ of the matrix $K$ approximate the eigenfunctions of the Laplace-Beltrami operator on the sampled manifold and thus can be used to parameterize the manifold~\cite{coifman2006diffusion}. Since $K$ is row-stochastic, the first eigenvector $\psi_0$ is trivial: all ones. The subsequent eigenvectors are called \emph{diffusion coordinates} and have corresponding eigenvalues $\lambda_1,\lambda_2,\cdots,\lambda_{m-1}$. The original data points are mapped to their diffusion coordinates as
	\begin{equation*}
		x\mapsto\big(\lambda_1^{\tau}\psi_1(x),\ldots,\lambda_{m-1}^{\tau}\psi_{m-1}(x)\big),
	\end{equation*}
	where $\psi_i(x)\in\mathbb{R}$ represents the entry of eigenvector $\psi_i$ corresponding to the point $x$ from the original data. In the following sections we take $\tau=0$. The distance between two mapped points is called {\em diffusion distance}, and it represents the similarity between two points in the original space. If two points are nearby in the diffusion space using a Euclidean metric, it implies that there are multiple short paths to transition from one point to the other in the original space. Similarly to PCA, if there is a spectral gap we can map our original data set to $k$ ``important" eigenvectors. However, attention must be paid to the fact that some eigenvectors may be higher harmonics of previous discovered ones~\cite{dsilva2018parsimonious}. The nonlinear manifold is parameterized by the few eigenvectors that correspond to the largest eigenvalues and that are not themselves such higher harmonics. These diffusion coordinates are the important intrinsic variables that parameterize the nonlinear manifold and indicate its dimensionality.
	
\begin{figure}[p]
		\centering
		\begin{subfigure}[]{0.67\linewidth}
			\includegraphics[width=\linewidth]{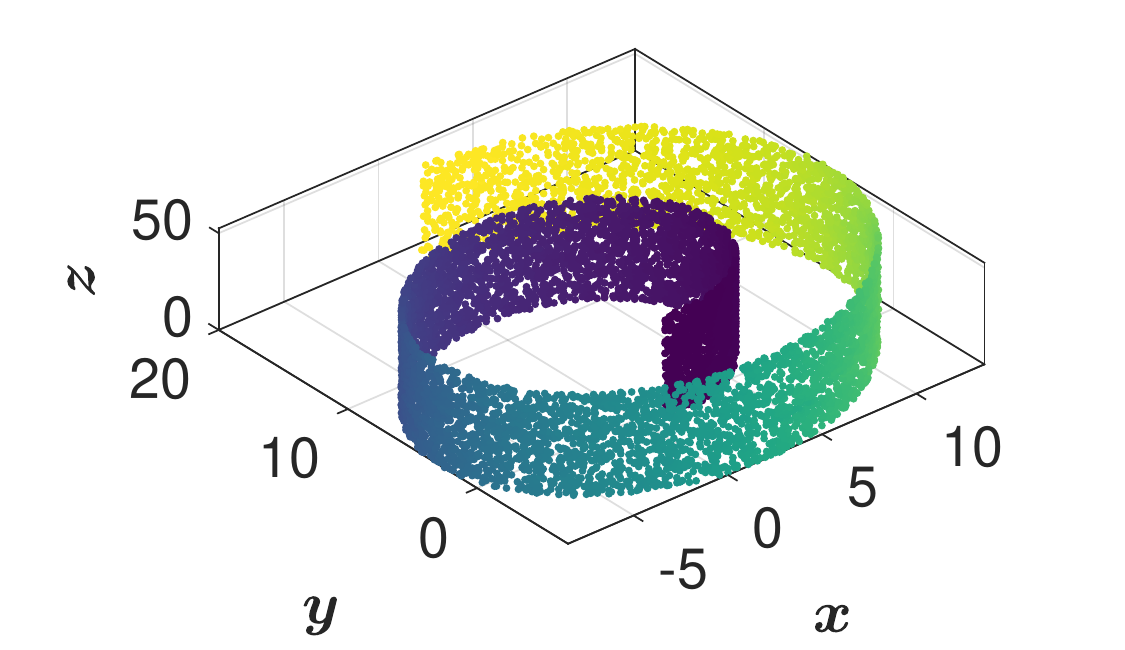}
			\caption{The first (nontrivial) diffusion coordinate $\psi_1$ parameterizes arclength along the roll.}
			\label{fig:swissrollpsi1}
		\end{subfigure}
		\qquad
		\begin{subfigure}[]{0.67\linewidth}
			\includegraphics[width=\linewidth]{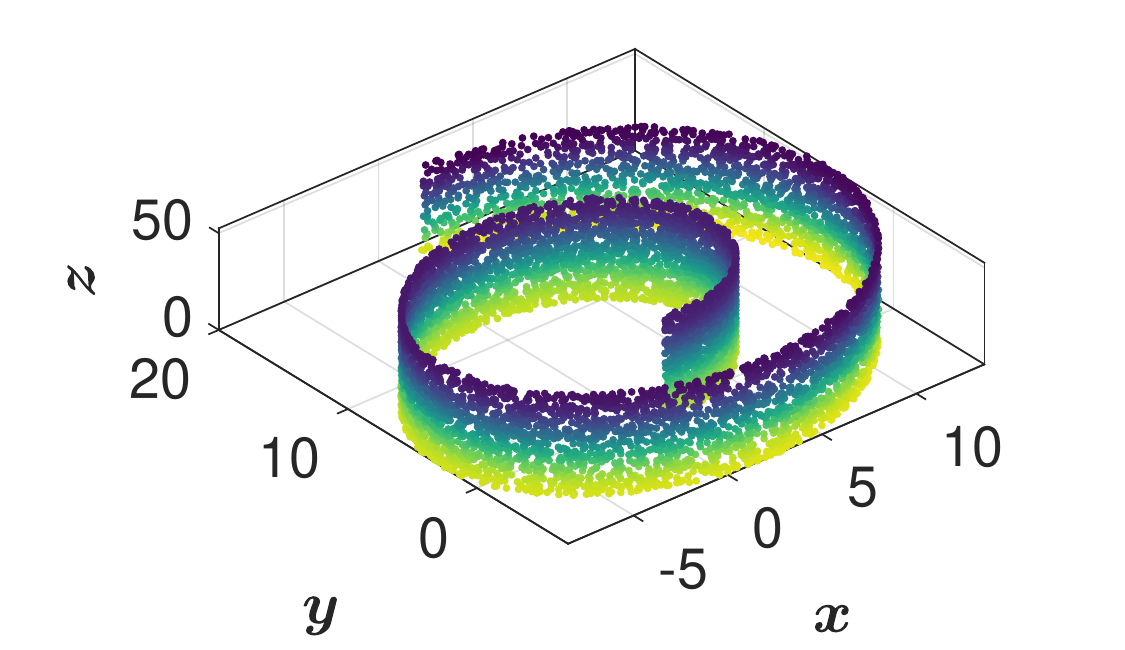}
			\caption{The second diffusion coordinate $\psi_2$ parameterizes height.}
			\label{fig:swissrollpsi2}
		\end{subfigure}
		\qquad
		\begin{subfigure}[]{0.67\linewidth}
			\includegraphics[width=\linewidth]{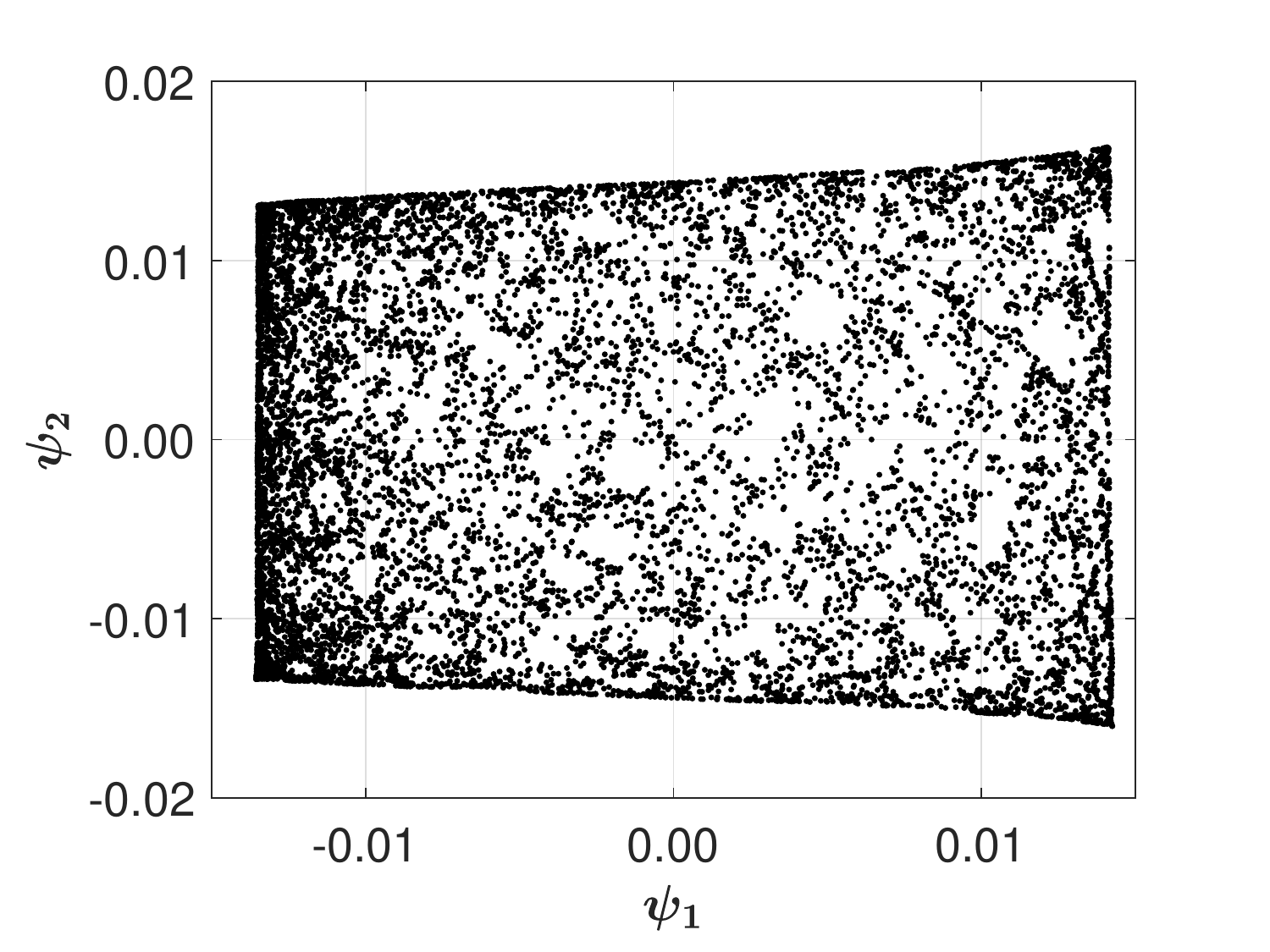}
			\caption{The two-dimensional nature of the data set is uncovered in the diffusion map space.}
			\label{fig:swissrollunrolled}
		\end{subfigure}
		\caption{Applying Diffusion Maps to the Swiss roll data set.}
		\label{fig:swissroll}
	\end{figure}

	We illustrate DMaps by applying it to a ``Swiss roll" data set. This is a three-dimensional data set, but only two variables are sufficient to describe every point: the height and the arclength along the roll. Figure~\ref{fig:swissrollpsi1} shows the data set colored by the first non-trivial eigenvector that parameterizes the arclength, while Figure~\ref{fig:swissrollpsi2} shows the data set colored by the second non-trivial eigenvector that parameterizes the height. Figure~\ref{fig:swissrollunrolled} shows the data set ``unrolled'' in the diffusion map space.
	
	In the case that the original data set comes from independent stochastic processes, as in our time series simulator, but we observe it through some nonlinear transformation $y = f(x)$, we can retrieve the original manifold using Mahalanobis distances~\cite{singer2008non,dsilva2013nonlinear}. It can be shown that if $x_i,x_j$ are two data points in the original space and $y_i,y_j$ are their nonlinear transformations, then
	\begin{equation*}
		\lVert x_i-x_j\rVert^2=\frac{1}{2}(y_i-y_j)^\top\left[(JJ^\top)^{-1}(y_i)+(JJ^\top)^{-1}(y_j)\right](y_i-y_j)+ O\big(\lVert y_i-y_j\rVert^4\big),
	\end{equation*}
	where $J$ is the Jacobian matrix of the transformation. In practice, the matrix $JJ^\top$ is approximated by a covariance matrix which is estimated by running several short bursts of our simulator around each data point. Since the manifold has a lower-dimensionality than the observed space, the covariance matrix will be rank deficient and {\em a pseudo-inverse must be used} to compute the Mahalanobis distances.
	
	\subsection{Parameter inference}
	We mentioned above that time series from the SA algorithm can be considered as corresponding to those of an effective stochastic differential equation. Hence, an essential component of our algorithm is the estimation of drift and diffusion coefficients of a stochastic process from local path data. Assume the one-dimensional stochastic process $dx(t)=h(x(t))\,dt+\sigma(x(t))\,dW$, with $W$ a standard Brownian motion. The coefficients can be estimated either from their statistical definitions~\cite{gradivsek2000analysis}, \emph{i.e.},
	\begin{equation}\label{eq:statdefn}
		\begin{split}
			&h\big(x(t)\big)=\lim_{\tau\to0}\frac{\big\langle x(t+\tau)-x(t)\big\rangle}{\tau}, \\
			&\sigma^2\big(x(t)\big)=\lim_{\tau\to0} \frac{\big\langle\big(x(t+\tau)-x(t)\big)^2\big\rangle}{\tau}, \\
		\end{split}
	\end{equation}
	or using the Generalized Method of Moments (GMM)~\cite{hansen1982large,chan1992empirical}, where moment conditions can be easily derived from an Euler-Maruyama discretization of the stochastic process. The above methods are more fitting if the stochastic process is realized as multiple short trajectories starting from the same initial conditions. On the other hand, if we are given a single, long trajectory then {\em maximum likelihood methods}~\cite{ait2002maximum,ait2008closed} are more suitable.
	
	The maximum likelihood estimator (MLE) for the drift coefficients is known to be asymptotically unbiased~\cite{Pavl2014}, \emph{i.e.}, as the length of the observed path increases, the MLE converges to the true values of the coefficients that appear in the drift. This is no longer true in the presence of a multiscale structure, \emph{i.e.}, when we want to estimate parameters in a stochastic coarse-grained model, given observations of the slow variable from the full dynamics. Indeed, it was shown rigorously in~\cite{papav-2007,PavlSt06} that in this case the MLE becomes asymptotically biased and that subsampling at an appropriate rate, between the two characteristic time scales of the dynamics, is needed in order to estimate accurately the parameters in the coarse-grained model. This is particularly relevant for us, since the optimization/estimation methodology is based on the assumption of the existence of a reduced model that describes accurately the system we are interested in. See, for example, the ODE driven by a sped-up Lorenz ’63 ODE that is studied in Section 3.4.
	
	\section{Results \& discussion}
	\subsection{An illustration: One-dimensional ``effective'' optimization}
	Before delving into the complete algorithm, which involves estimation of effective gradients in the low-dimensional embedding, we proceed with a simpler example where the objective function is two-dimensional
	but the optimization process can be effectively one-dimensional. Consider an objective function given by
	\begin{equation}\label{eq:bayesmodel}
		f(x,y)=C\exp\left[-25\left(x^2+y^2\right)^2+216(x^2+y^2)-0.05\sqrt{(x-2)^2+y^2}\,\right]\,,
	\end{equation}
	where $C$ is a proportionality constant selected to normalize the maximum objective value to unity. Since~(\ref{eq:bayesmodel}) arises as the posterior density of a Bayesian parameter estimation, our algorithm in reality minimizes $-f$, but we present the results from the perspective of maximizing $f$.
	
\begin{figure}[htb]
		\centering
		\includegraphics[width=\linewidth]{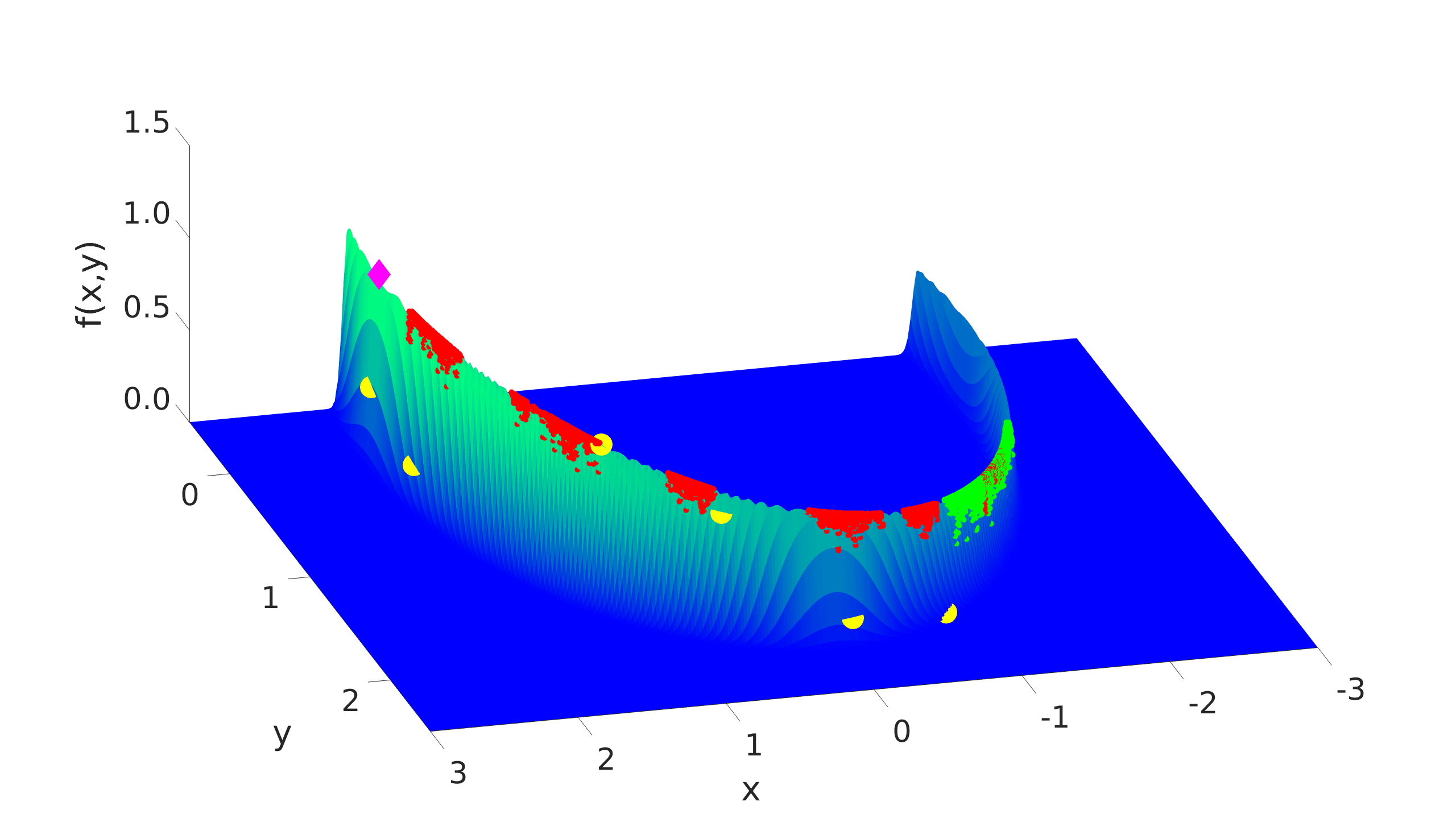}
		\caption[Coarse-grained optimization of a Bayesian model]{One complete run of the algorithm that approaches the global maximum. Six total ``coarse iterations" (shown in red) were performed. In addition, a single run of SA using the same number of function evaluations is shown in green. Our algorithm visibly approaches the maximum much faster.}
		\label{fig:bayesianExample}
	\end{figure}

	Plotting the objective (cf. Figure \ref{fig:bayesianExample}) reveals that $f\approx0$ on the entire plane except for a thin region near the circle $x^2+y^2=4$, where the function exhibits a sharply peaked ``ridge'' of greater function values. If we attempt to maximize $f$ via Simulated Annealing (or, at constant temperature, RWMH) the algorithm will quickly be attracted to the ridge and then slowly proceed along it towards the maximum. Trial points far away from the ridge are usually worse than the current accepted point and are highly likely to be rejected. Hence, the accepted points lie on an essentially one-dimensional curve in the parameter space.
	
	We can apply Diffusion Maps to this data, and the first diffusion coordinate will parameterize the curve. Once we have obtained the one-dimensional embedding, we can extrapolate in the direction that the objective function increases. The last step involves returning from the diffusion space back to the original parameter space. This procedure is called ``lifting", and a variety of methods are available, such as Laplacian Pyramids~\cite{dsilva2013nonlinear}, Geometric Harmonics~\cite{coifman2006geometric,lafon2006data} and Radial Basis Functions~\cite{chiavazzo2014reduced}. We will use geometric harmonics in the present work. After obtaining the projected point back in parameter space, we perform another short run of SA from that point and repeat the procedure. To summarize the algorithm:
	\begin{enumerate}
		\item Pick an initial point, possibly near the ``ridge" of the objective function.
		\item Run a short burst of Simulated Annealing until a prescribed number of points has been accepted (1000 here).
		\item Discard any outliers that may be far away from the ridge.
		\item Apply Diffusion Maps to the remaining two-dimensional data set and obtain a nonlinear embedding. The embedding is one-dimensional and the diffusion coordinate can be thought of as corresponding to the arclength along the ridge.
		\item In the diffusion space, project to a new point that is in the direction that increases the objective function.
		\item Lift the new point back to full space via geometric harmonics. The resulting point is expected to be close to the ridge, but even if it is not, the new SA run will quickly be attracted to it.
		\item Return to Step 2 and repeat the procedure.
	\end{enumerate}
	
\begin{figure}[htb]
		\centering
		\includegraphics[width=\linewidth]{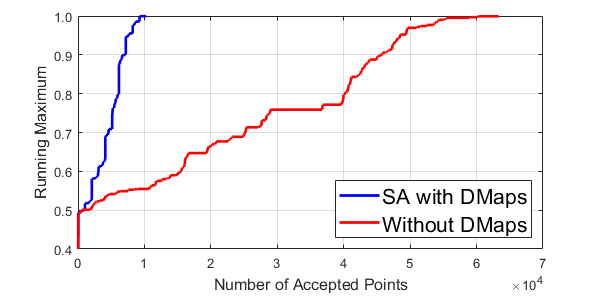}
		\caption{A comparison of the evolving maximum objective value for both methods. SA combined with DMaps needs only a fraction of the total function evaluations compared to simple SA.}
		\label{fig:bayesianComp}
	\end{figure}

	The results for this illustrative example are shown in Figure~\ref{fig:bayesianExample}. Each short burst of SA is shown as a ``cloud" of red points. The lifted points are shown in yellow, and we can see that geometric harmonics perform well in this case, as the lifted points are still close to the ridge. The maximum of the function is depicted by the purple diamond. Additionally, a single run of SA using the same total number of function evaluations was performed (in green). Figure \ref{fig:bayesianComp} compares the running maximum objective value achieved by the two approaches. It is clear that, for this simple illustration, SA/RWMH combined with Diffusion Maps approaches the maximum substatially more quickly. Building on the body of ideas developed in \cite{gear2010computing,gear2003projective} as well as \cite{vanden2003fast,weinan2007heterogeneous}, the combination of an ``inner optimizer" with data mining of its local results, followed by taking larger ``outer" optimization steps in the identified reduced space, has the potential to significantly accelerate the overall computational optimization.
	
	\subsection{Coefficient estimation of the effective SDE}
	In this section we will demonstrate how one can estimate the theoretically expected drift and diffusion coefficients after the trajectories of a stochastic process have been
	transformed using DMaps. As we mentioned before, these coefficients correspond to an effective gradient and will provide us with an approximation of the ``correct"
	ascent (resp. descent) direction along which to optimize (maximize, resp. minimize) in the low-dimensional space.
	We begin with a two-dimensional SDE, analogous to the Langevin equation:
	\begin{equation}
		\begin{split}\label{eq:langevin_analogue}
			dx & =\mu(x,y)\,dt+\sqrt{2T}\,dW_1 \\
			dy & =\nu(x,y)\,dt+\sqrt{2T\,}dW_2,
		\end{split}
	\end{equation}
	where $[W_1\,\,W_2]^T = \mathbf{W}$ are independent Brownian motions. Assume now that the data in the original space $x,y$ are transformed by being observed through the leading Diffusion Map coordinates, and the trajectories are now written in terms of these diffusion coordinates $\psi_1(x,y)$ and $\psi_2(x,y)$. In order to rewrite our system of SDEs in terms of the new variables, we apply the multidimensional It\^{o}'s lemma~\cite{oksendal2003stochastic,Pavl2014}, \textit{e.g.}, for $\psi_1$ we have:
	\small
	\begin{equation*}
		\begin{split}
			d\psi_1 & =\left((\nabla\psi_1)^\top\left[\begin{array}{c}\mu \\ \nu\end{array}\right]+\frac{1}{2}\text{Tr}\left[\Sigma^\top(\mathbf{H}\psi_1)\Sigma\right]\right)dt+(\nabla\psi_1)^\top\Sigma\,d\mathbf{W}\\
			& = \left[\left(\frac{\partial\psi_1}{\partial x}\mu+\frac{\partial\psi_1}{\partial y}\nu\right)+T\left(\frac{\partial^2\psi_1}{\partial x^2}+\frac{\partial^2\psi_1}{\partial y^2}\right)\right]dt \\
			& \qquad+\sqrt{2T}\sqrt{\left(\frac{\partial\psi_1}{\partial x}\right)^2+\left(\frac{\partial\psi_1}{\partial y}\right)^2}\,d\tilde{W_1},
		\end{split}
	\end{equation*}
	\normalsize
	where $\Sigma=\sqrt{2T}\,I$ is the covariance matrix from~(\ref{eq:langevin_analogue}), $\mathbf{H}\psi_1$ is the Hessian matrix of second partial derivatives of $\psi_1$, and $\tilde{W_1}$ is a new Brownian motion.
	
	In order to simplify the estimation, we set up a two-dimensional grid in the $x,y$ space. The partial derivatives that are required for the theoretical computation of coefficients are approximated numerically at the grid points using centered differences. From every grid point, $N$ trajectories are simulated via SA for a specified time $\Delta t$, which is also the time step in the estimation computed via~(\ref{eq:statdefn}). The simulation time step between successive points on a trajectory is $\delta t=0.1\,\Delta t$.
	
	The data set that is then ``passed" to the DMaps algorithm, in order to find the new embedding, consists of the cloud of final points from each trajectory, the initial grid points, and all points where the partial derivatives are estimated. Afterwards, we have a new grid in the $\psi_1,\psi_2$ space, with every partial derivative estimated on this grid. We can assume that in a small neighborhood of every grid point the partial derivatives of the diffusion coordinates are approximately constant. Given that, we perform a separate estimation at each grid point of the following system of SDEs:
	\begin{equation*}
		\begin{split}
			d\psi_1 & =\theta_1\,dt+\theta_3\,d\tilde{W_1}\\
			d\psi_2 & =\theta_2\,dt+\theta_4\,d\tilde{W_2},
		\end{split}
	\end{equation*}
	where $\theta_1$ and $\theta_3$ correspond to
	\begin{equation}\label{eq:theoretical}
		\begin{split}
			\theta_1 & =\left[\left(\frac{\partial\psi_1}{\partial x}\mu+\frac{\partial\psi_1}{\partial y}\nu\right)+T\left(\frac{\partial^2\psi_1}{\partial x^2}+\frac{\partial^2\psi_1}{\partial y^2}\right)\right] \\
			\theta_3 & =\sqrt{2T}\sqrt{\left(\frac{\partial\psi_1}{\partial x}\right)^2+\left(\frac{\partial\psi_1}{\partial y}\right)^2}\,,
		\end{split}
	\end{equation}
	and similarly for $\theta_2$ and $\theta_4$. Thus, we obtain values for each $\theta_i,\,i=1,2,3,4$ at every grid point and fit a polynomial along the grid, \textit{e.g.}, using a quadratic fit:
	\begin{equation*}
		\theta_1(\psi_1,\psi_2) \approx p_{00} + p_{10}\psi_1 + p_{01}\psi_2 + p_{20}\psi_1^2 + p_{11}\psi_1\psi_2 + p_{02}\psi_2^2.
	\end{equation*}
	If the grid is local (small) enough and the bursts are contained within the grid for the most part, we can assume that the objective function could be approximated locally by a linear surface which has a constant gradient, \textit{i.e.}:
	\begin{equation*}
		\begin{split}
			&\mu(x,y) = \mu_0 \\
			&\nu(x,y) = \nu_0.
		\end{split}
	\end{equation*}
	
	For our illustrative example we will use the function $f(x,y)=x+2y$, which has a constant gradient, on an orthogonal grid along the coordinate axes. We use an $8\times10$ grid with limits $[0,1.5]\times[0,1.2]$, and we use $N=150$ trajectories at every grid point, each run for $\Delta t=0.01$. The time step of the simulation is $\delta t=10^{-3}$.
	
\begin{figure}[p]
		\centering
		\begin{subfigure}[]{0.45\linewidth}
			\includegraphics[width=\linewidth]{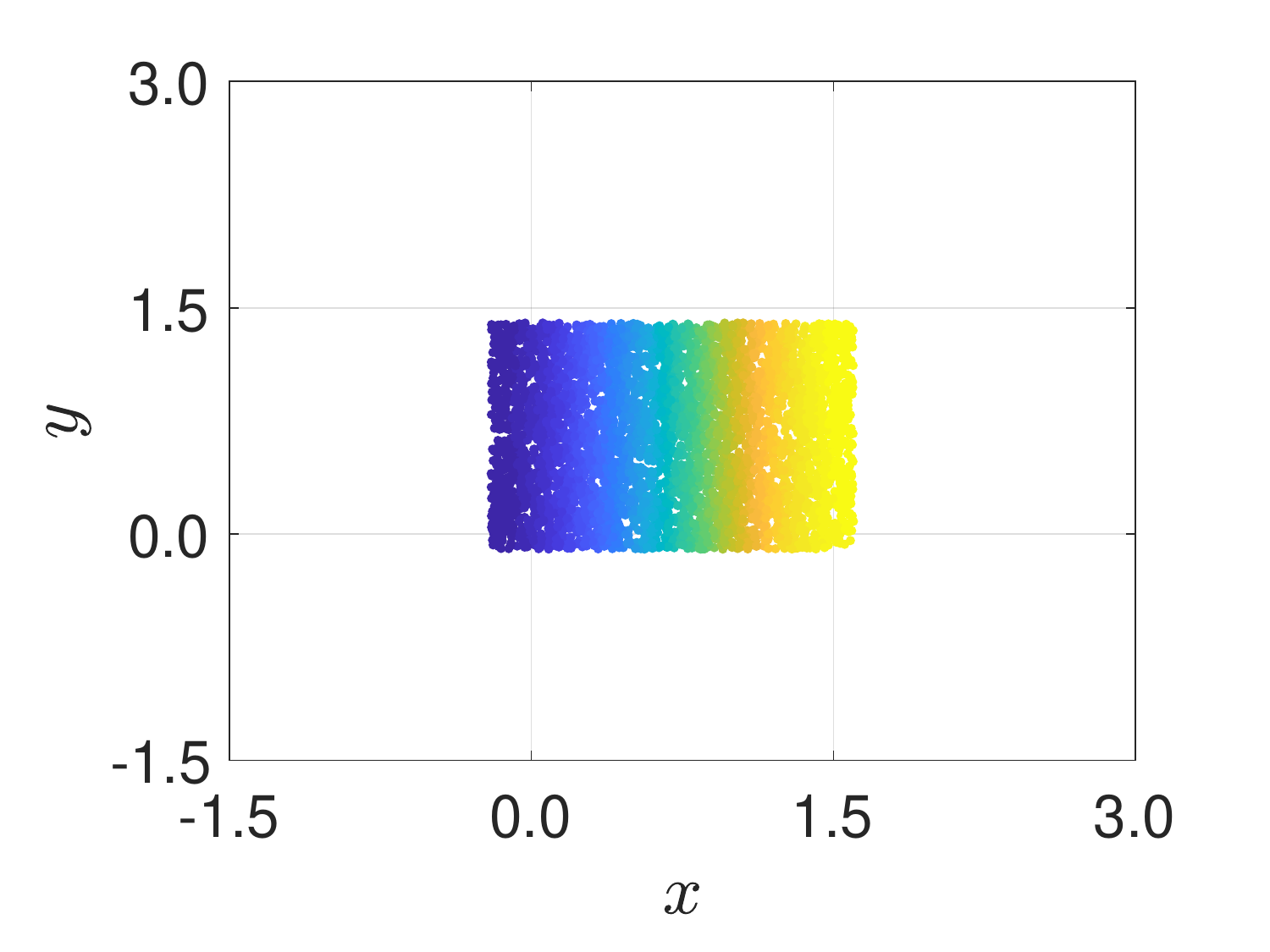}
			\label{fig:2Dpsi1}
		\end{subfigure}
		\qquad
		\begin{subfigure}[]{0.45\linewidth}
			\includegraphics[width=\linewidth]{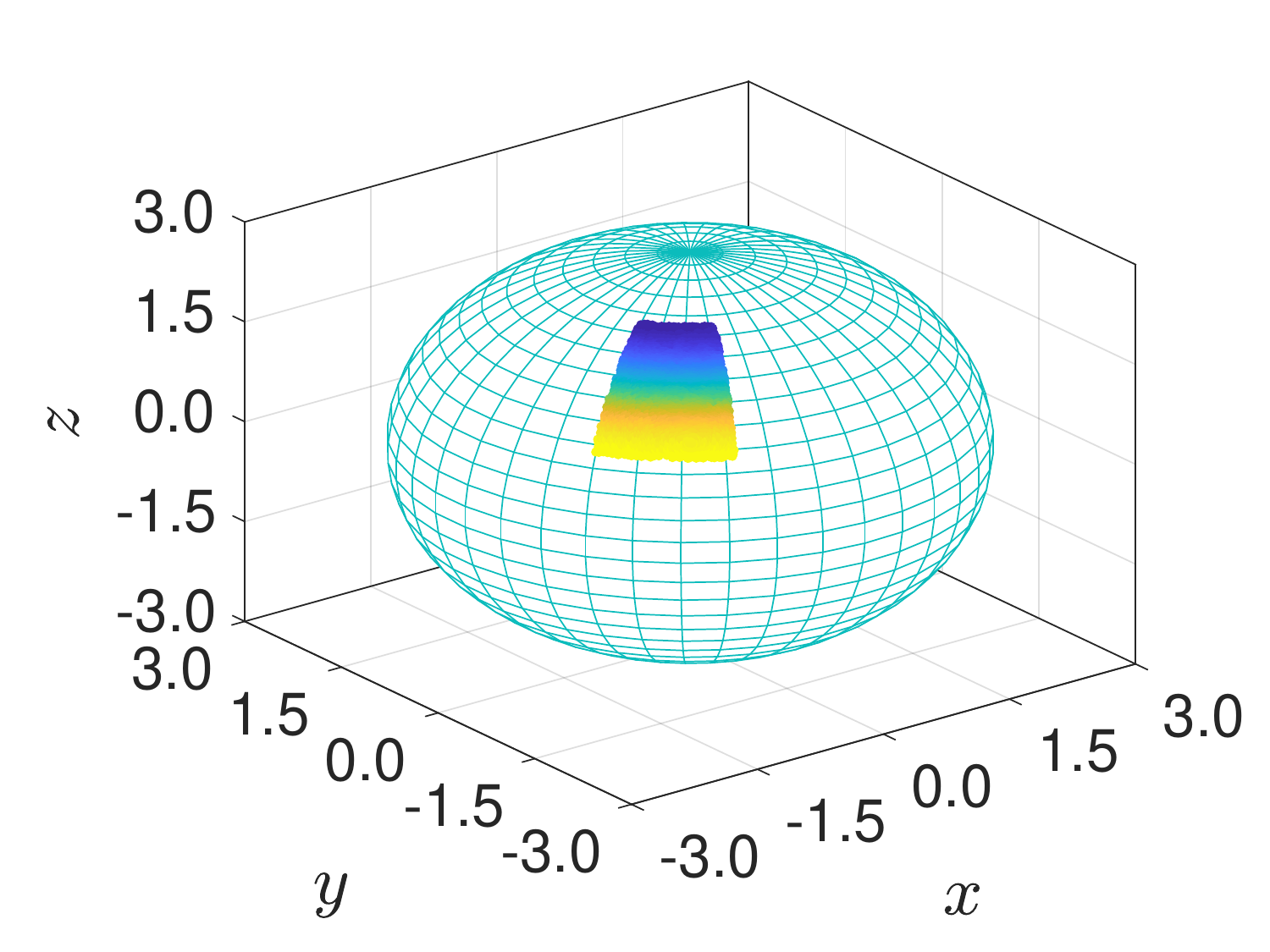}
			\label{fig:3Dpsi1}
		\end{subfigure}
		\qquad
		\begin{subfigure}[]{0.45\linewidth}
			\includegraphics[width=\linewidth]{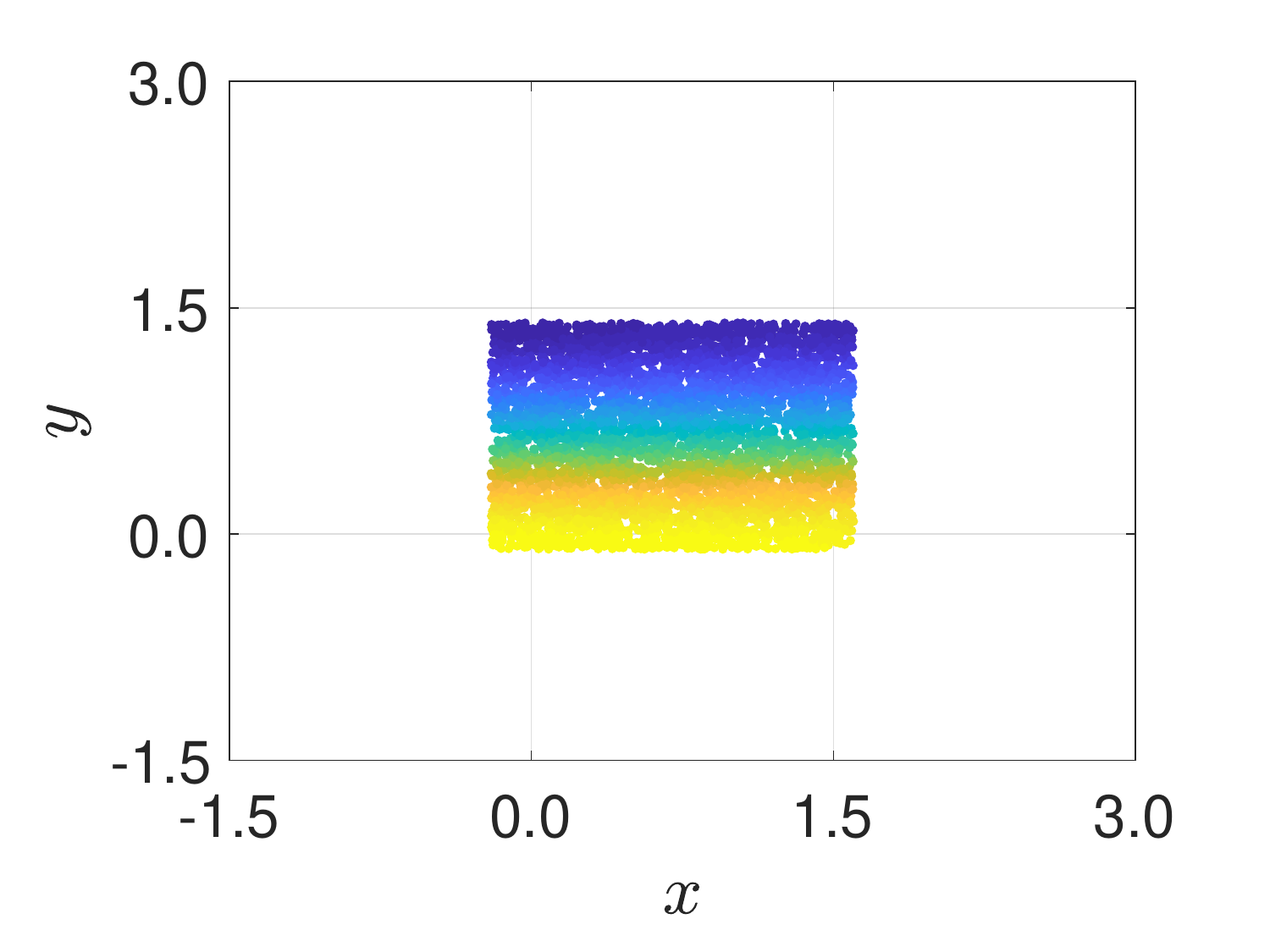}
			\label{fig:2Dpsi2}
		\end{subfigure}
		\qquad
		\begin{subfigure}[]{0.45\linewidth}
			\includegraphics[width=\linewidth]{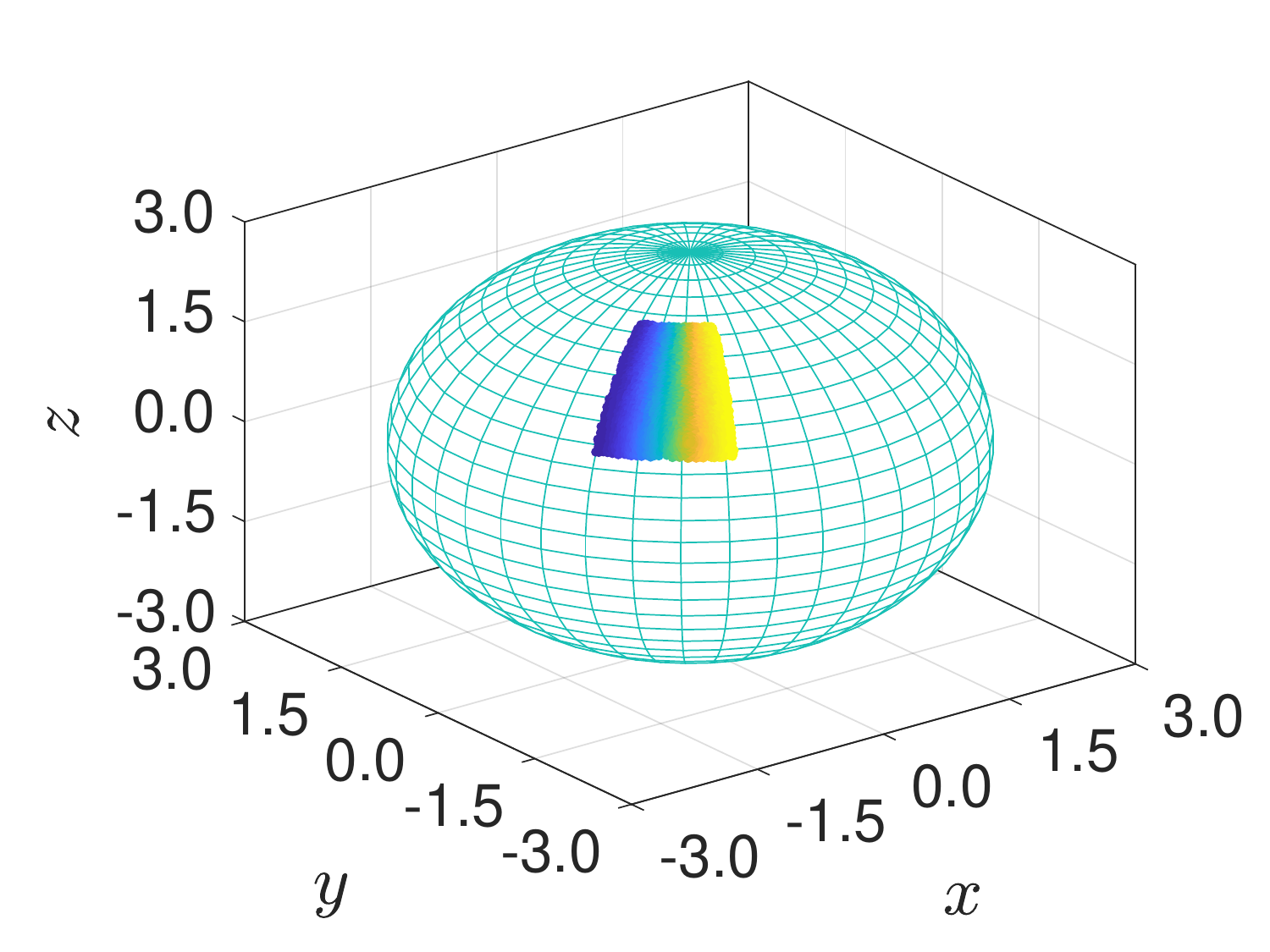}
			\label{fig:3Dpsi2}
		\end{subfigure}
		\qquad
		\begin{subfigure}[]{0.45\linewidth}
			\includegraphics[width=\linewidth]{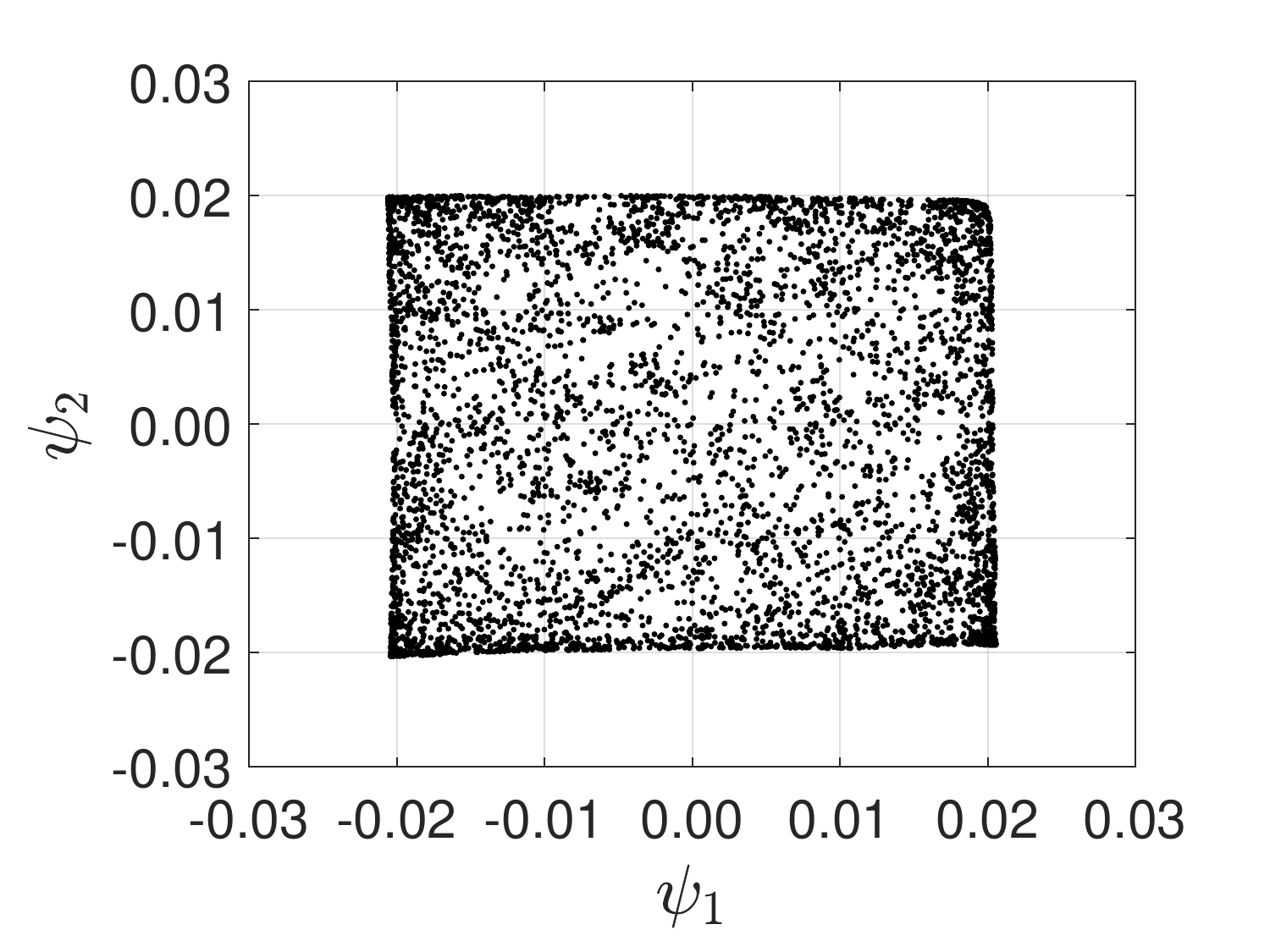}
			\caption{}
			\label{fig:2Dembedding}
		\end{subfigure}
		\qquad
		\begin{subfigure}[]{0.45\linewidth}
			\includegraphics[width=\linewidth]{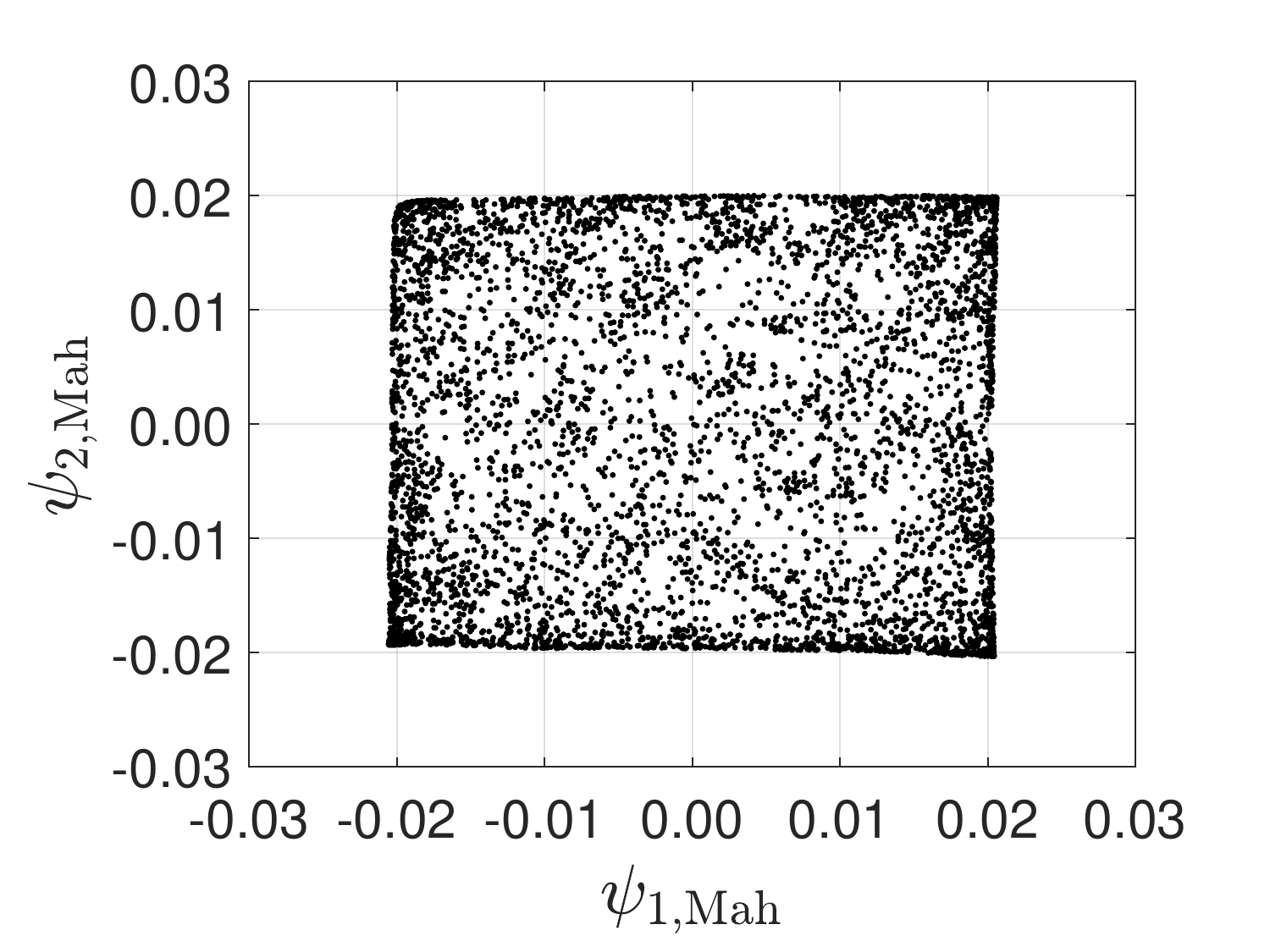}
			\caption{}
			\label{fig:3Dembedding}
		\end{subfigure}
		\caption[Diffusion Maps applied to bursts on a rectangular strip]{Diffusion maps applied to (A) the original 2D data with Euclidean distances, and (B) the transformed data with Mahalanobis distances. In both cases, the inherent dimensionality is recovered in the first two eigenvectors, as the coloring of the data by the leading Diffusion Map coordinates corresponding to these two eigenvectors show.}
		\label{fig:estimation}
	\end{figure}

	To illustrate that the estimation is accurate even if we observe the process through a nonlinear transformation, we transform a region of the $(x,y)$ space that contains the trajectories and map it onto a spherical surface, obtaining a new $(x,y,z)$ space.
	We apply Diffusion Maps with Euclidean distances to the original data set and Diffusion Maps with Mahalanobis distances to the transformed data set. Figure~\ref{fig:2Dembedding} shows the original data set colored by the two diffusion coordinates and the new embedding. Of course, since the original data set is two-dimensional, no dimensionality reduction is achieved in this case. Figure~\ref{fig:3Dembedding} shows the transformed data set colored by the diffusion coordinates. In this case, Mahalanobis distances enable us to retrieve the original, orthogonal, two-dimensional embedding.
	
\begin{figure}[htb]
		\centering
		\begin{subfigure}[]{0.8\linewidth}
			\includegraphics[width=\linewidth]{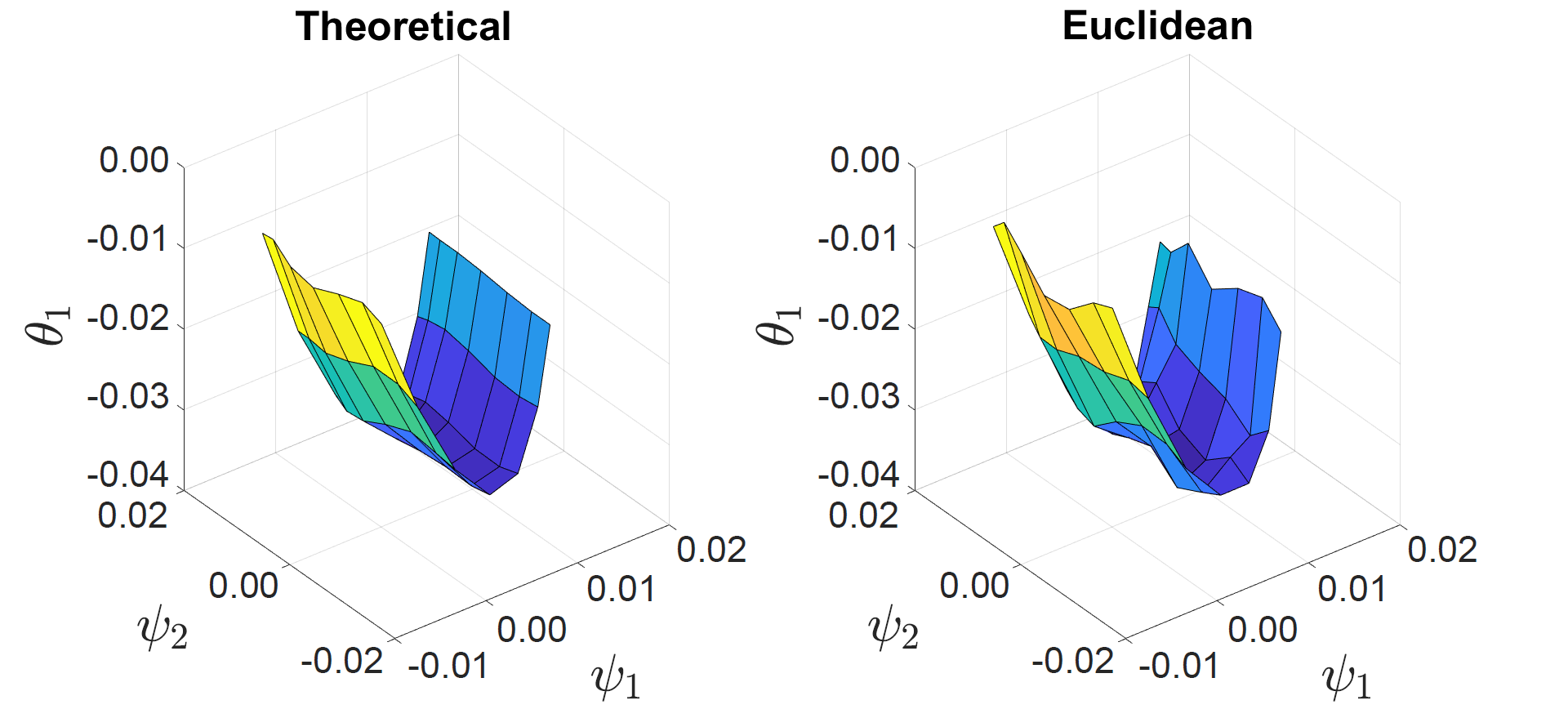}
			\label{fig:theta1_coeffA}
		\end{subfigure}
		\begin{subfigure}[]{0.8\linewidth}
			\includegraphics[width=\linewidth]{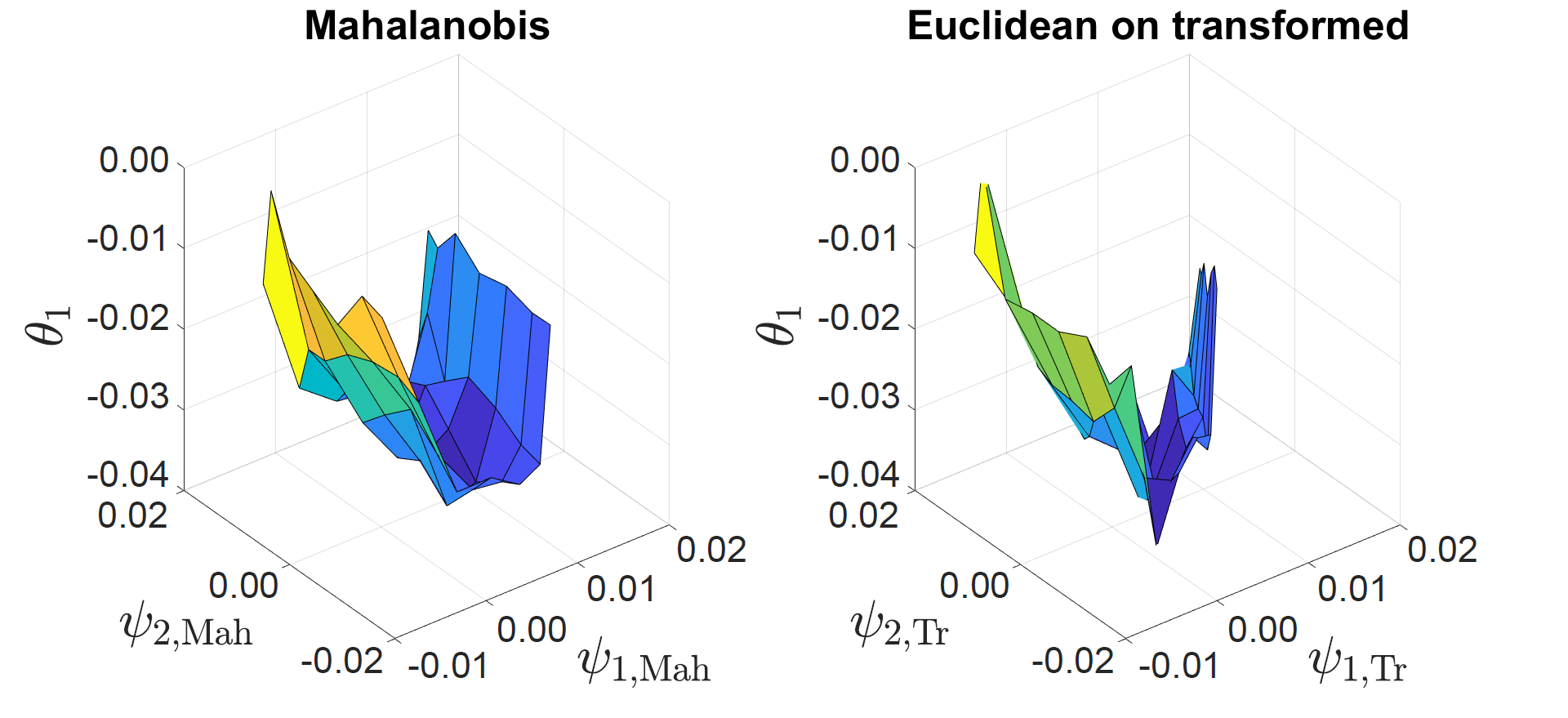}
			\label{fig:theta1_coeffB}
		\end{subfigure}
		\caption[Estimation of first drift coefficient]{Estimation of the first drift coefficient $\theta_1$. ``Theoretical'' are obtained numerically from~(\ref{eq:theoretical}), ``Euclidean'' via DMaps on the original data, and ``Mahalanobis'' from DMaps on the transformed data. The last subplot shows results from Euclidean DMaps on transformed data, which, as expected, yields incorrect estimates.}
		\label{fig:theta1}
	\end{figure}

\begin{figure}[htb]
		\centering
		\begin{subfigure}[]{0.8\linewidth}
			\includegraphics[width=\linewidth]{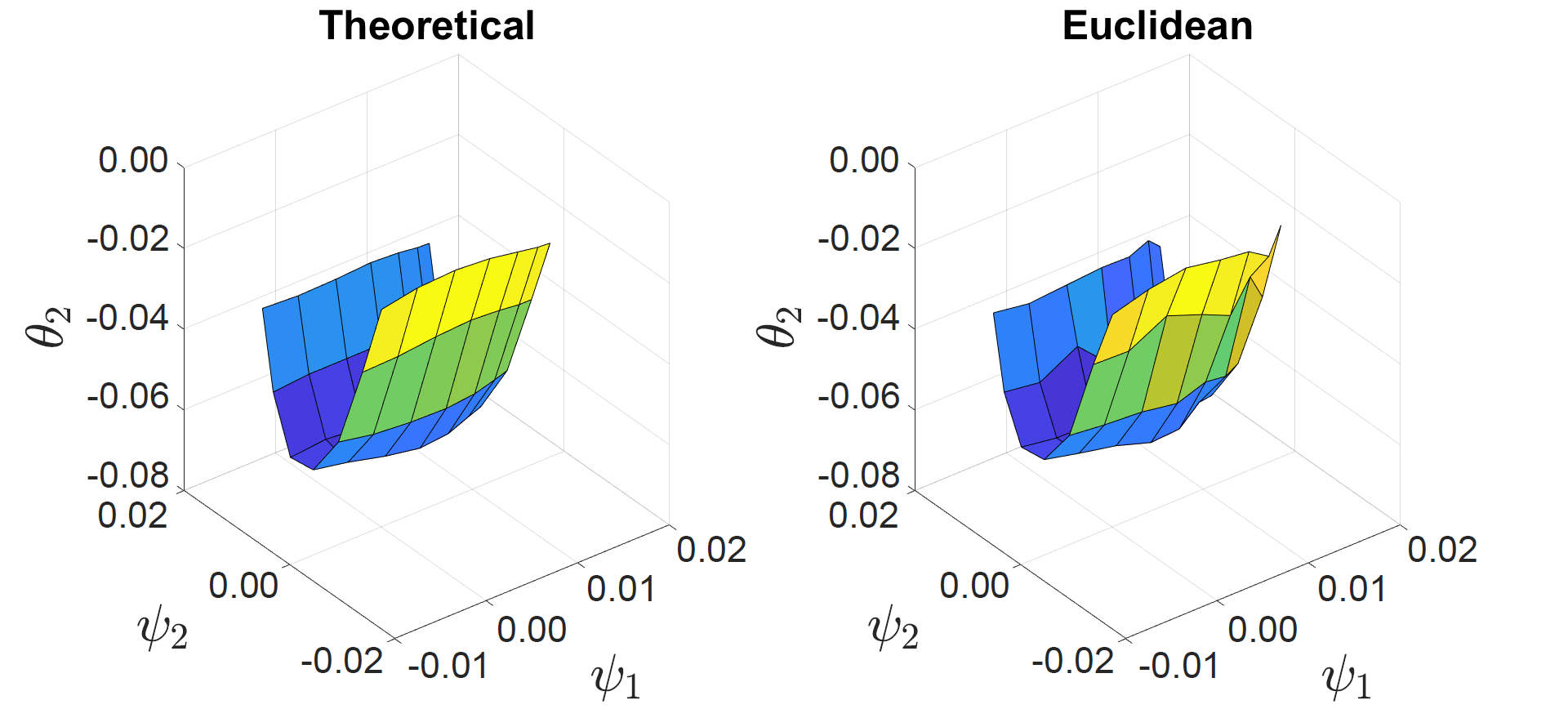}
			\label{fig:theta2_coeffA}
		\end{subfigure}
		\begin{subfigure}[]{0.8\linewidth}
			\includegraphics[width=\linewidth]{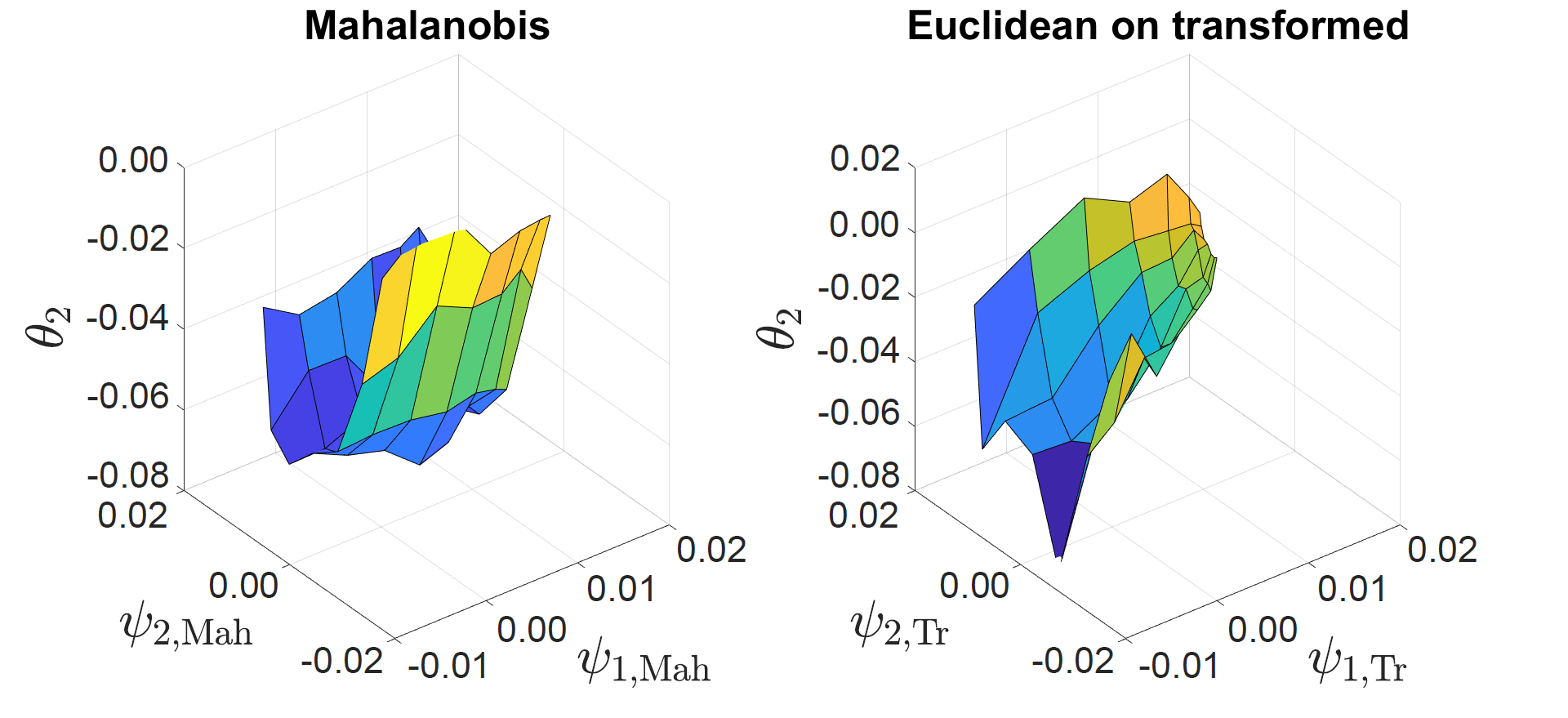}
			\label{fig:theta2_coeffB}
		\end{subfigure}
		\caption[Estimation of second drift coefficient]{Estimation of the second drift coefficient $\theta_2$. The subplots are analogous to those in Figure~\ref{fig:theta1}.}
		\label{fig:theta2}
	\end{figure}

	After estimation, some of the drifts and diffusivities along the grid approach zero. In order to avoid these degeneracies, we discard the points where this occurs and fit the coefficients $\theta_i$ to the rest of the grid. Figures~\ref{fig:theta1} and~\ref{fig:theta2} show the results for drift coefficients $\theta_1,\theta_2$. Similar results are obtained for diffusion coefficients $\theta_3,\theta_4$.
	
	\subsection{Optimization on a cylinder}
	Having shown that we can estimate the correct parameters in the Diffusion Map space whether we observe the original manifold or some invertible transformation of it, we now apply our algorithm to a three-dimensional objective function with an attracting, slow, two-dimensional manifold.
	
	In a cylindrical coordinate scheme (\textit{i.e.}, $x(r,\theta,z)=r\cos\theta,\;y(r,\theta,z)=r\sin\theta$, $z(r,\theta,z)=z$), we define:
	\begin{equation*}
		f(r,\theta,z)=\frac{k_1}{2}(r-R)^2+h(\theta)+\frac{k_2}{2}z^2,
	\end{equation*}
	where $k_1,\,k_2>0$ determine to what extent trajectories are attracted to the circle defined by the intersection of the plane $z=0$ with the cylinder having radius $R$ and axis $z$. We used parameter values $(k_1,\,k_2,\,R)=(10^4,\,20,\,5/\pi)$ and implemented $f$ programmatically so as to accept Cartesian coordinates and internally convert to cylindrical.
	The function
	\begin{equation*}
		h(\theta)=-1.2+3.4\cos^2\left(\theta\right)-0.59\cos\left(\theta\right)-1.1\sin\left(\theta\right)
	\end{equation*}
	determines an asymmetric double well potential with a local minimum close to $\theta = -\pi/2$ and a global minimum around $\theta = \pi/2$. Any trajectory away from the cylinder is quickly attracted to it, and then the search of the parameter space proceeds along the cylinder surface. The new algorithm builds on the previously presented one, but now also includes parameter estimation in Diffusion Map space.
	
	\vspace*{\baselineskip}
	\noindent
	\paragraph*{\underline{Algorithm}}
	\begin{enumerate}
		\item Initialize a local grid around a starting point and simulate ensembles of short trajectories starting at every grid point.
		\item Apply DMaps to the data set to obtain the low-dimensional embedding.
		\item Estimate SDE coefficients on the grid points.
		\item Fit a polynomial to the coefficients along the grid.
		\item Integrate  the system of ODEs $d\psi_1=\theta_1(\psi_1,\,\psi_2)\,dt,\,d\psi_2=\theta_2(\psi_1,\,\psi_2)\,dt$ forward in time. This is analogous to using a gradient descent algorithm.
		\item Lift the resulting point to full space and use it as a new starting point.
		\item Repeat until the estimated coefficients in Diffusion Map space approach zero.
	\end{enumerate}
	
\begin{figure}[p]
		\centering
		\begin{subfigure}[]{0.8\linewidth}
			\includegraphics[width=\linewidth]{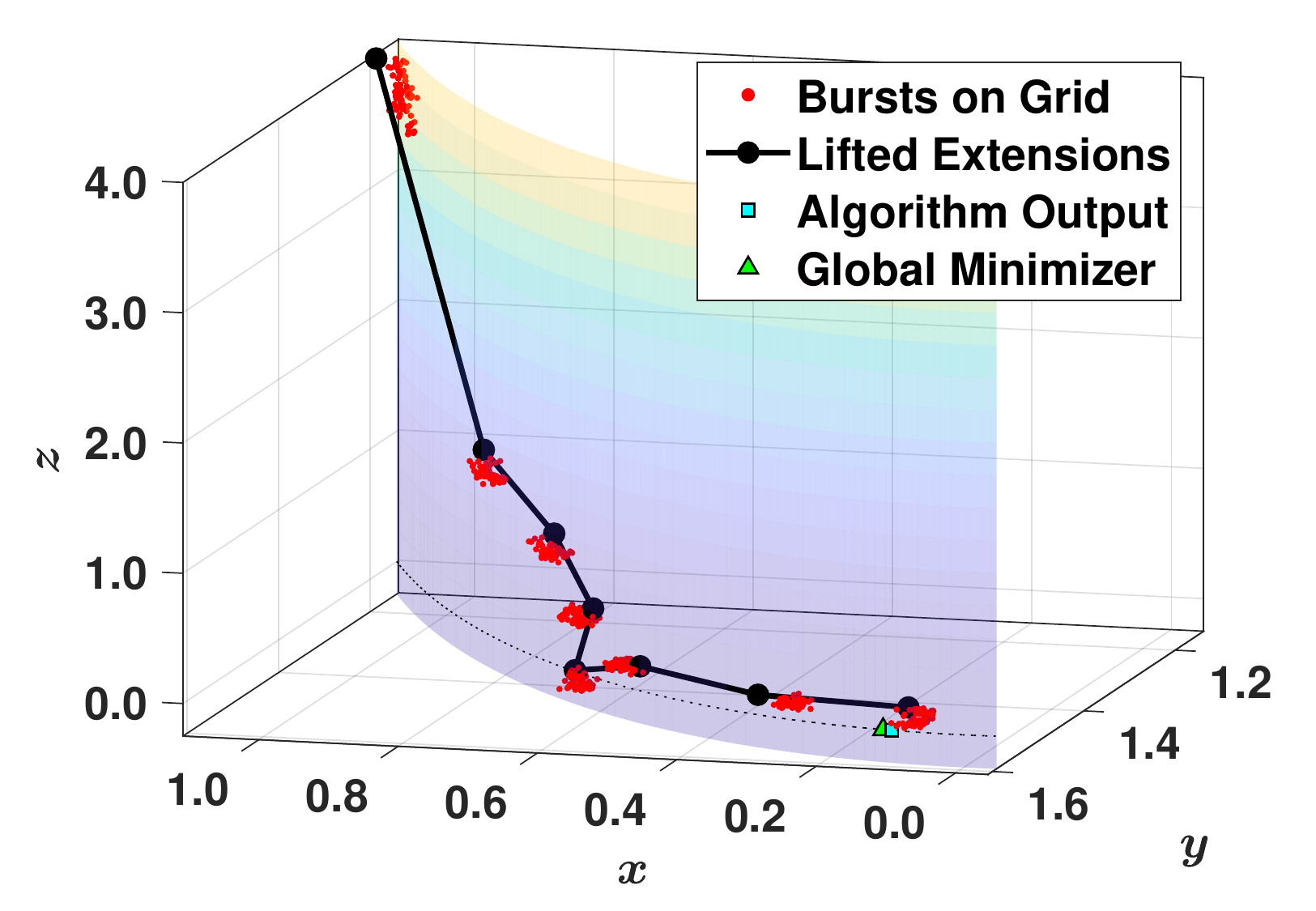}
			\label{fig:runInit}
		\end{subfigure}
		\begin{subfigure}[]{0.8\linewidth}
			\includegraphics[width=\linewidth]{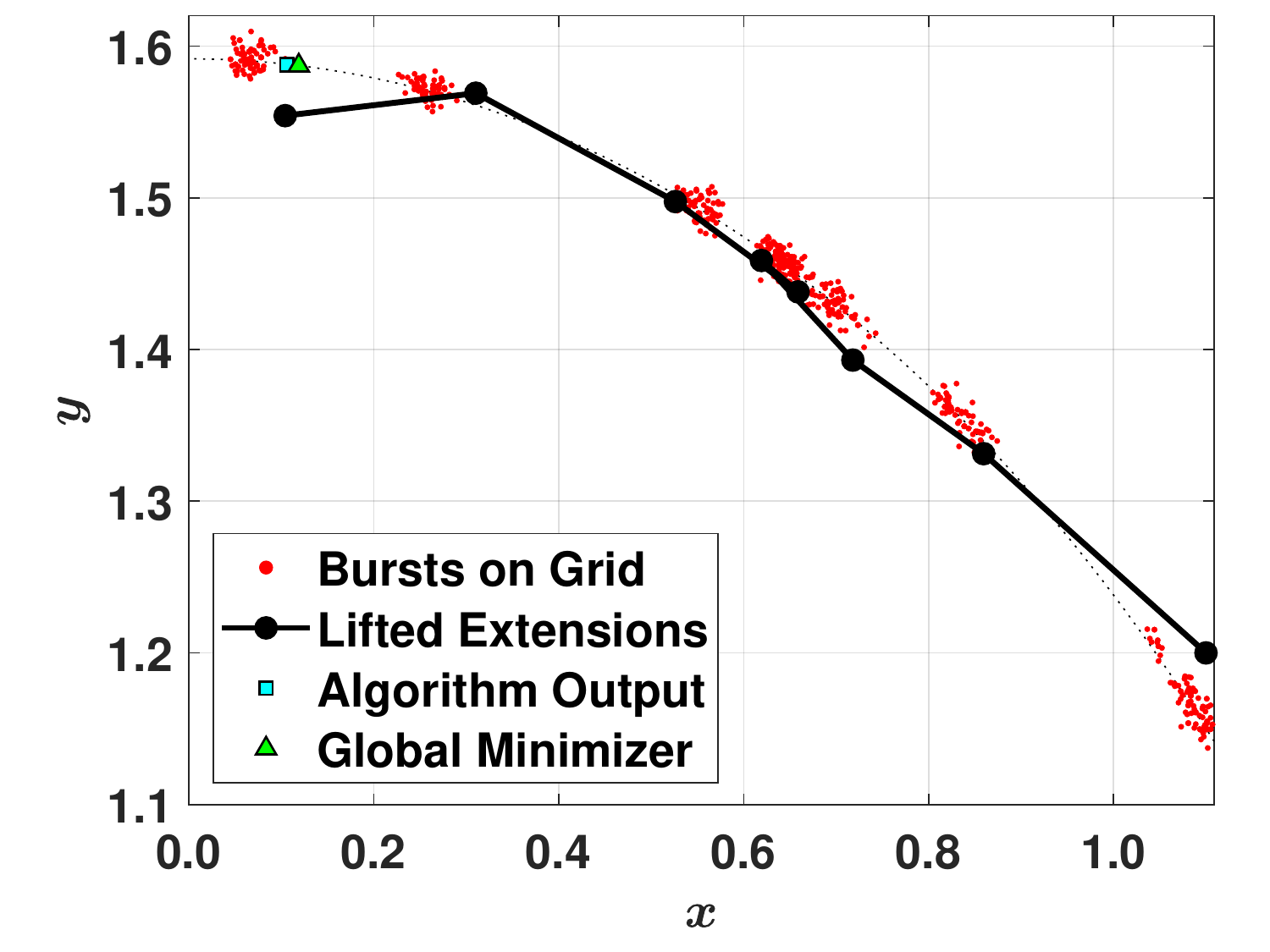}
			\label{fig:runFin}
		\end{subfigure}
		\caption[Coarse-grained optimization of a function embedded on a cylindrical surface]{Coarse-grained (two-dimensional) optimization on a cylindrical surface in three dimensions.  One complete run of the algorithm, illustrating the short bursts of trajectories and new points after being lifted by geometric harmonics. Observe that, although the lifting does not perfectly locate the cylinder, the burst trajectories are quickly attracted back to it. The second plot is a top-down view of the first.}
		\label{fig:2runcylinder}
	\end{figure}

	\vspace*{\baselineskip}
	One must of course pay attention to the length of the integration step, to avoid extrapolating too far away from the grid in 3D space, since the coefficients are known there only locally. In general, we observed that linear fits perform better and follow the negative gradient direction at larger distances from the grid. Figure~\ref{fig:2runcylinder} shows two snapshots of the algorithm run as it approaches the minimum near $\theta=\pi/2$.
	
	It is important to note that the above procedure is a simple first attempt at a general outline of how to perform optimization in the reduced space. Many possible improvements can be made to reduce redundant computations and make the algorithm more efficient for practical applications. The first issue that should be addressed is how many short runs of the optimizer are required to obtain sufficient information about an effective gradient. The grid setup in the previous example is potentially redundant, and estimation at only a few points may be sufficient to obtain a good ascent/descent direction. Another question is how far along the manifold one can usefully project. A line search method~\cite{kelley1999iterative} could be implemented here, though one should keep in mind that the effective gradient is estimated locally and that, the farther away we extend along the manifold, the less accurate the lifting procedure becomes. We are systematically exploring these considerations for future publication.
	
	\subsection{Fast chaotic noise}
	Our estimation procedure can also be applied in cases where the underlying stochasticity of the system is not due to a Wiener process but arises from deterministic chaos. Consider an ODE driven by one of the components of the Lorenz system:
	\begin{equation}
		\begin{aligned}\label{eq:noiseSDE}
			\frac{dx}{dt} & =A(x-x^3)+\frac{\lambda}{\ep}y_2,\\
			\frac{dy_1}{dt} & =\frac{10}{\ep^2}(y_2-y_1), \\
			\frac{dy_2}{dt} & =\frac{1}{\ep^2}(28y_1-y_2-y_1y_3),\\
			\frac{dy_3}{dt} & =\frac{1}{\ep^2}(y_1y_2-\frac{8}{3}y_3).
		\end{aligned}
	\end{equation}
	
	It can be shown that the approximate dynamics for the slow variable are given by the following SDE~\cite{KKP2015,melbourne2011note}:
	\begin{equation*}
		dx=A(x-x^3)\,dt+\sqrt{\sigma}\,dW.
	\end{equation*}
	In the following simulations we use parameter values of $A=1$, $\lambda = 2/45$, and $\ep = \sqrt{0.001}$. Assuming Diffusion Maps yields a diffusion coordinate $\psi(x)$ that is one-to-one with the slow variable $x$, we can write an SDE for it using It\^{o}'s Lemma:
	\begin{equation}\label{eq:ito4phi}
		d\psi=\left(A(x-x^3)\frac{d\psi}{dx}+\frac{\sigma}{2}\frac{d^2\psi}{dx^2}\right)dt+\sqrt{\sigma}\frac{d\psi}{dx}dW.
	\end{equation}
	
	The simulation is set up as follows (cf.~\cite{krumscheid2013semiparametric,KPPK2015}): initially, the system is run long enough from initial conditions $(1,1,1,1)^\top$ so that it converges onto the Lorenz attractor ($t\approx0.1$). The end point of this initial simulation of the Lorenz system $(\mathbf{y}_0)$ will be used subsequently as our starting point for our short bursts. The starting points for $x$ are taken as equally-spaced points in the range $[-1.5,\;1.5]$; we chose to use 20 such points. In order to achieve faster separation of trajectories starting close to the same $x$ value, we perturb the starting point for each short burst. The actual initial conditions are given as
	\begin{equation*}
		\begin{split}
			x_{ic} & =x_0+0.01\frac{x_{\text{spacing}}}{2}z, \\
			\mathbf{y}_{ic} & =\mathbf{y}_0+\mathbf{z},
		\end{split}
	\end{equation*}
	where $z$ are standard normal variables and $x_{\text{spacing}}$ is the distance between starting points for $x$ before the perturbation.
	
\begin{figure}[htb]
		\centering
		\includegraphics[width=\linewidth]{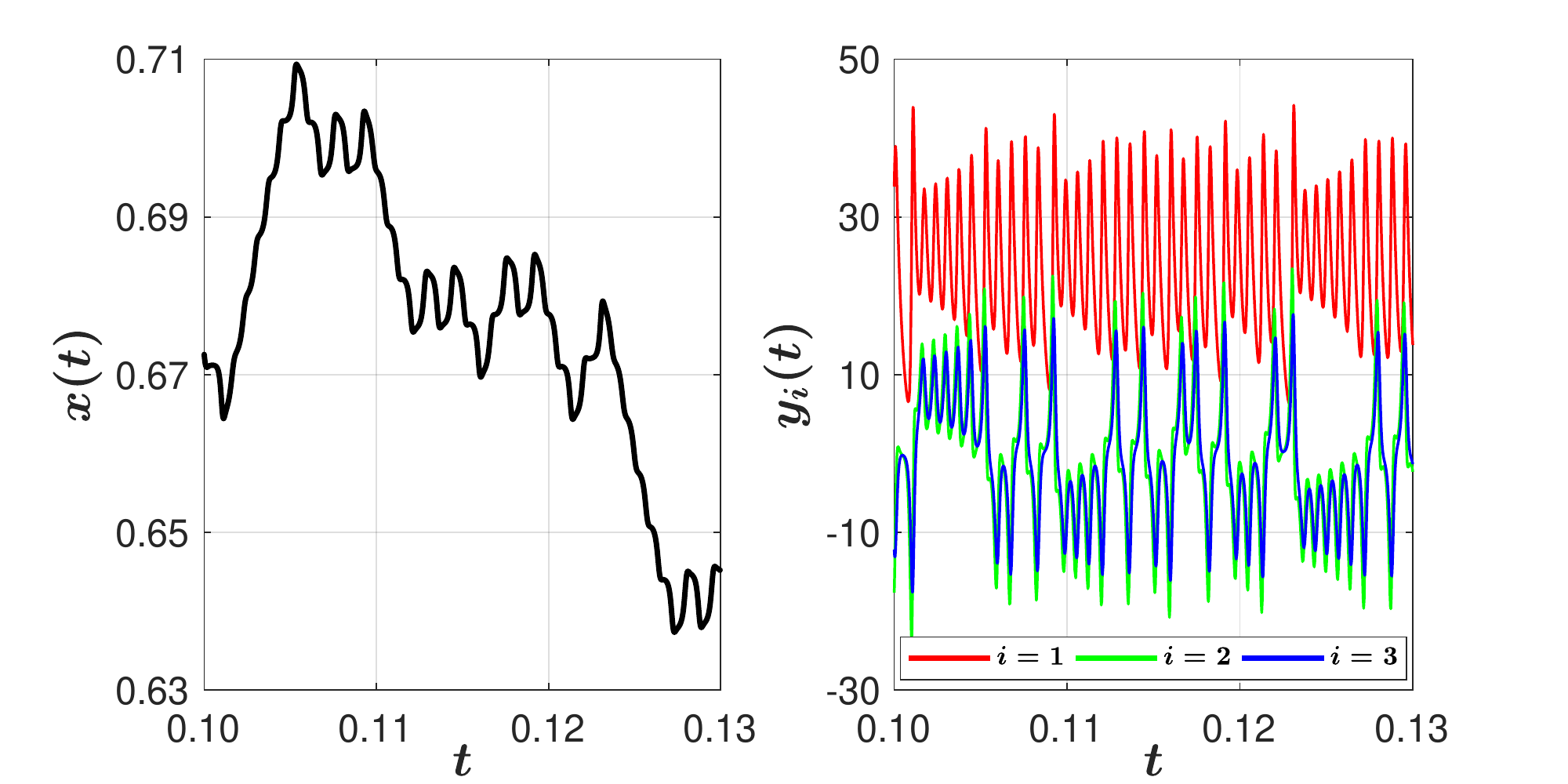}
		\caption{Simulation of a short trajectory initialized at perturbed initial conditions $(x_{ic},\mathbf{y}_{ic})$.}
		\label{fig:timescale}
	\end{figure}

	At each of our 20 starting points, 500 short trajectories are simulated. The system is integrated with time step $\delta t=10^{-3}$ for a duration $\Delta t=0.03$. Figure \ref{fig:timescale} shows the trajectories of both the fast and slow variables for one such burst. We see that the duration $\Delta t=0.03$ lies between the timescales of the fast (Lorenz) and slow (effective $x$) dynamics, so we avoid biasing our estimators for the coarse-grained model parameters.
	
	As before, we assume that the drift and diffusion coefficients are approximately constant at each starting point and estimate the following SDE at each starting point via GMM:
	\begin{equation*}
		dx=\theta_1\,dt+\theta_2\,dW.
	\end{equation*}
	Afterwards, we fit a polynomial of an appropriate degree to the estimated coefficients and retrieve the coefficients $\hat{A} \approx 0.9534$ and $\hat{\sigma} \approx 0.117$, values that are close to the ones reported in~\cite{krumscheid2013semiparametric}. The entire data set consists of the starting points, end points of each short burst, and points that are used to compute the derivatives of the diffusion coordinate. In this particular case, its dimension is $\mathbb{R}^{10060\times4}$.
	
\begin{figure}[htb]
		\centering
		\includegraphics[width=\linewidth]{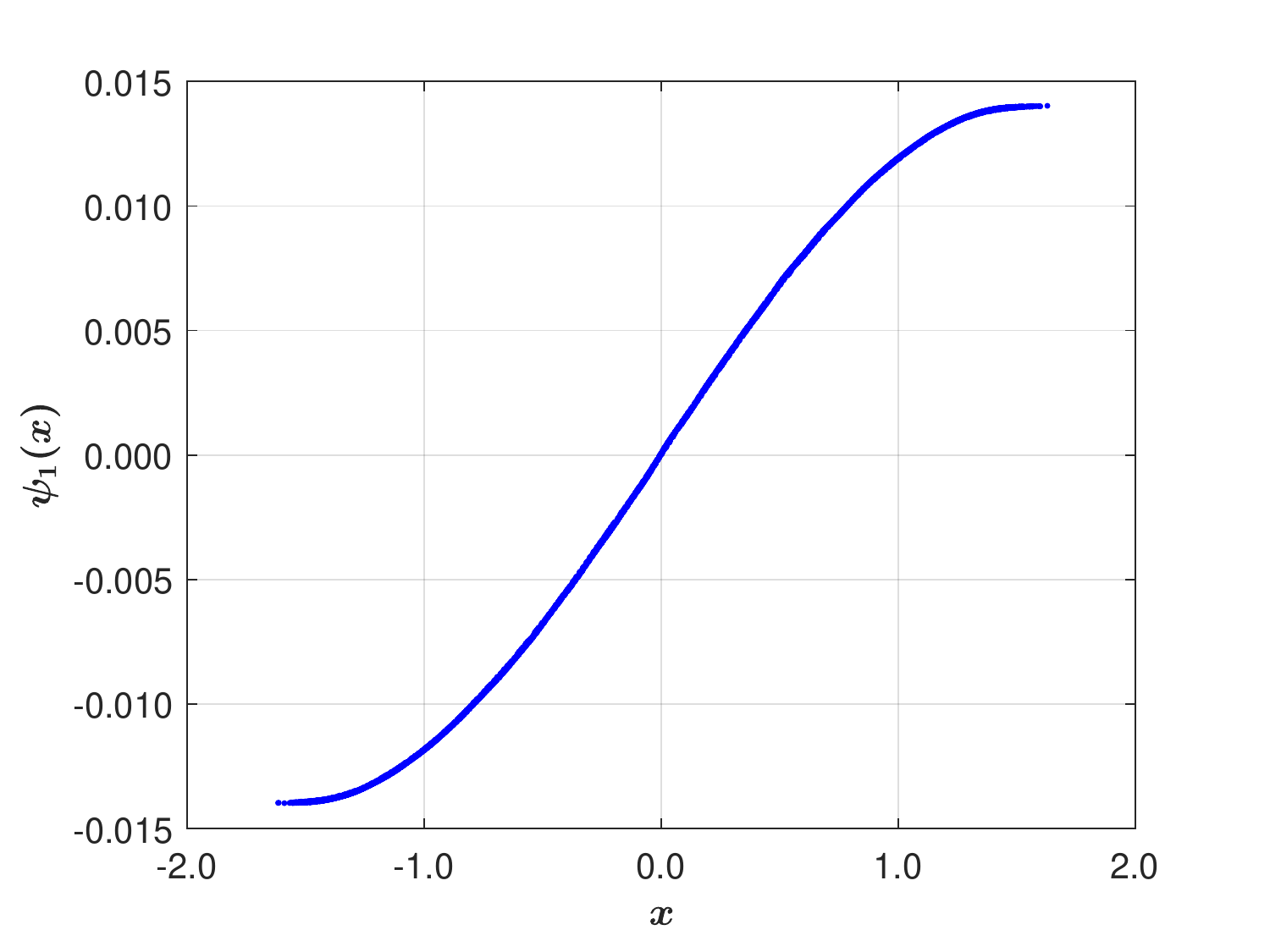}
		\caption{The first diffusion coordinate parameterizes $x$.}
		\label{fig:psi1}
	\end{figure}

	Applying regular DMaps to this data set yields a parameterization of the Lorenz attractor. In order to extract the slow variable, we apply DMaps using Mahalanobis distances. The local covariances of each data point are computed using 100 short simulations with duration $dt_{\text{cov}}=10^{-4}$. The pseudo-inverse of the covariance matrix is computed using a singular value decomposition (SVD):
	\begin{equation}
		C^\dagger=\sum_{m=1}^ds_m^{-1}v_mv_m^\top,\label{eq:svd}
	\end{equation}
	where $s$ are the singular values and $v$ the right-singular vectors. In this case, we use $d=n=4$, since we need the last singular value that corresponds to the slow variable to obtain the correct embedding. Using this setup, the first non-trivial diffusion coordinate parameterizes the slow variable $x$, as seen in Figure~\ref{fig:psi1}.
	
	While in theory the derivatives of the diffusion coordinate could be computed using central differences, the parameterization of $x$ is quite noisy and the estimated derivatives can be quite inaccurate. To ameliorate this difficulty, we used smoothing splines to fit a curve to the data and estimated the derivatives from the splines. The data and fit curve are shown in Figure~\ref{fig:psi1}.
	
\begin{figure}[p]
		\centering
		\includegraphics[width=\linewidth]{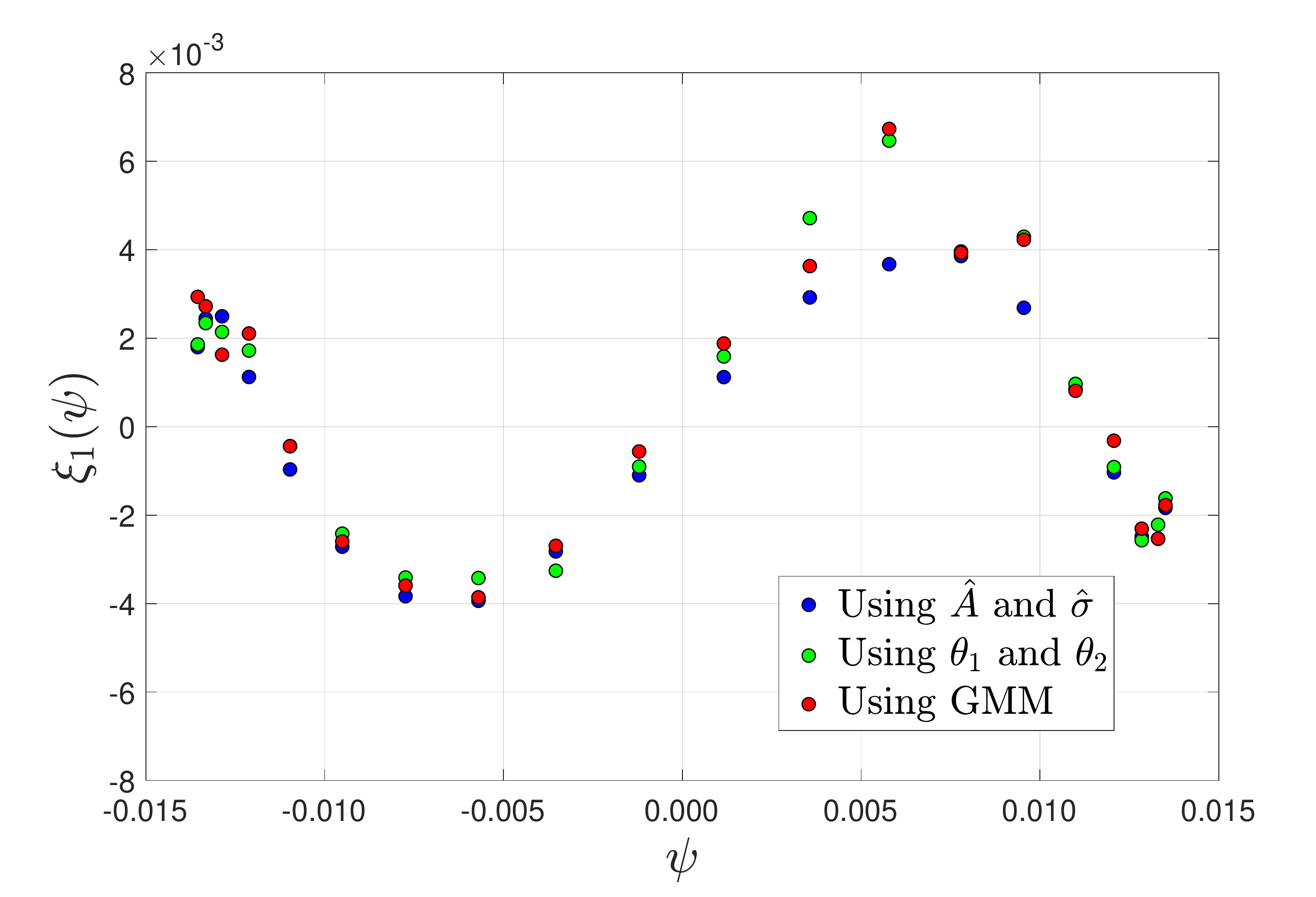}
		\caption{Estimation of the drift coefficient in diffusion map space.}
		\label{fig:ksi1}
	\end{figure}

\begin{figure}[p]
		\centering
		\includegraphics[width=\linewidth]{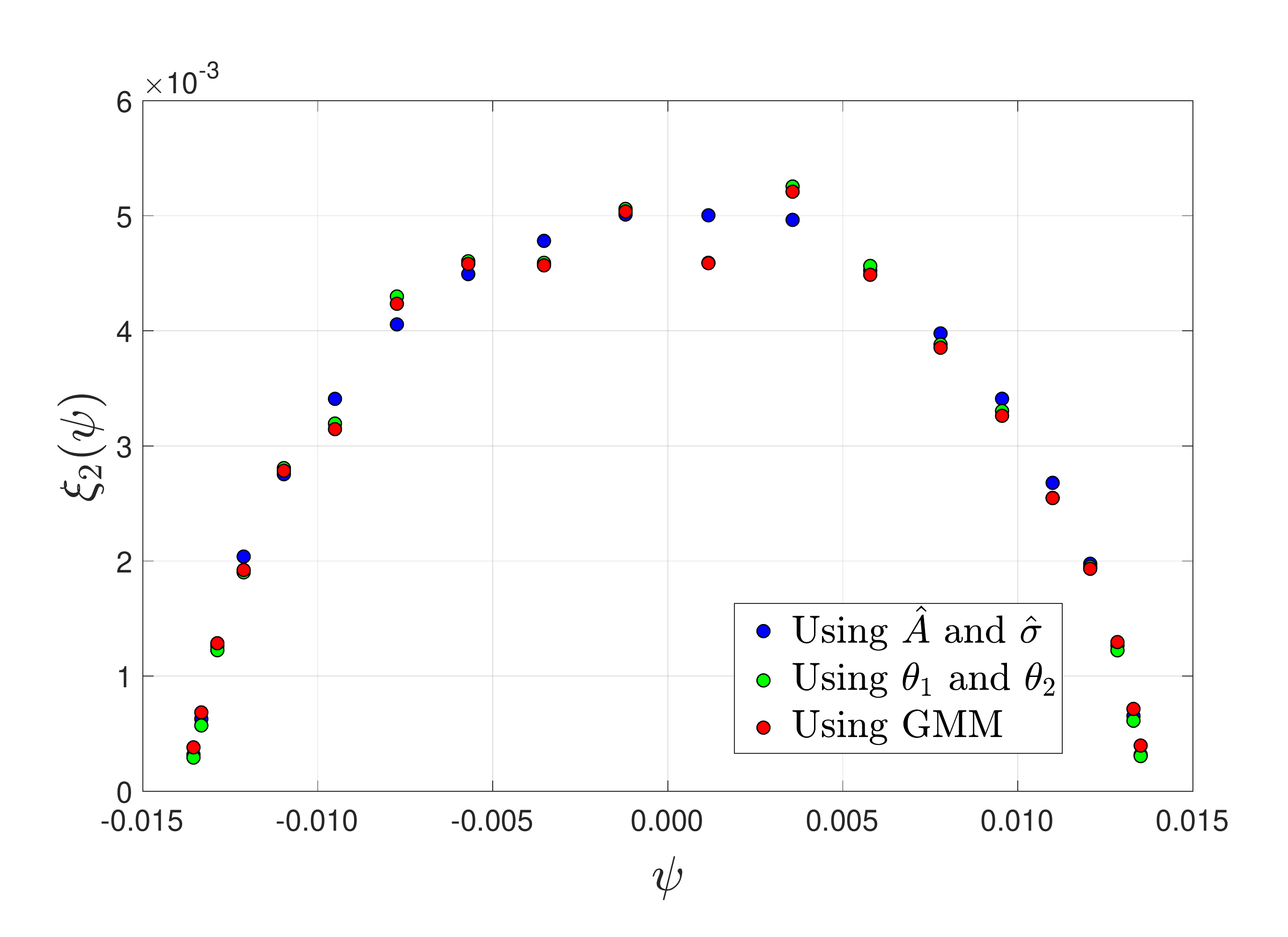}
		\caption{Estimation of the diffusion coefficient in diffusion map space.}
		\label{fig:ksi2}
	\end{figure}

	Using the new data set in Diffusion Map space we again assume constant local drift and diffusivity at each starting point and use GMM to fit the following SDE:
	\begin{equation*}
		d\psi=\xi_1\,dt+\xi_2\,dW.
	\end{equation*}
	This estimate is compared to (\ref{eq:ito4phi}) either using the estimates $\hat{A}$ and $\hat{\sigma}$:
	\begin{equation}
		d\psi=\left(\hat{A}(x-x^3)\frac{d\psi}{dx}+\frac{\hat{\sigma}}{2}\frac{d^2\psi}{dx^2}\right)dt+\sqrt{\hat{\sigma}}\,\frac{d\psi}{dx}\,dW,\label{eq:1}
	\end{equation}
	or directly using the estimates $\theta_1$ and $\theta_2$:
	\begin{equation*}
		d\psi=\left(\theta_1\frac{d\psi}{dx}+\frac{\theta_2^2}{2}\frac{d^2\psi}{dx^2}\right)dt+\theta_2\frac{d\psi}{dx}\,dW.
	\end{equation*}
	The results are shown in Figure~\ref{fig:ksi1} and Figure~\ref{fig:ksi2}. 
	
	
\begin{figure}[htb]
		\centering
		\includegraphics[width=\linewidth]{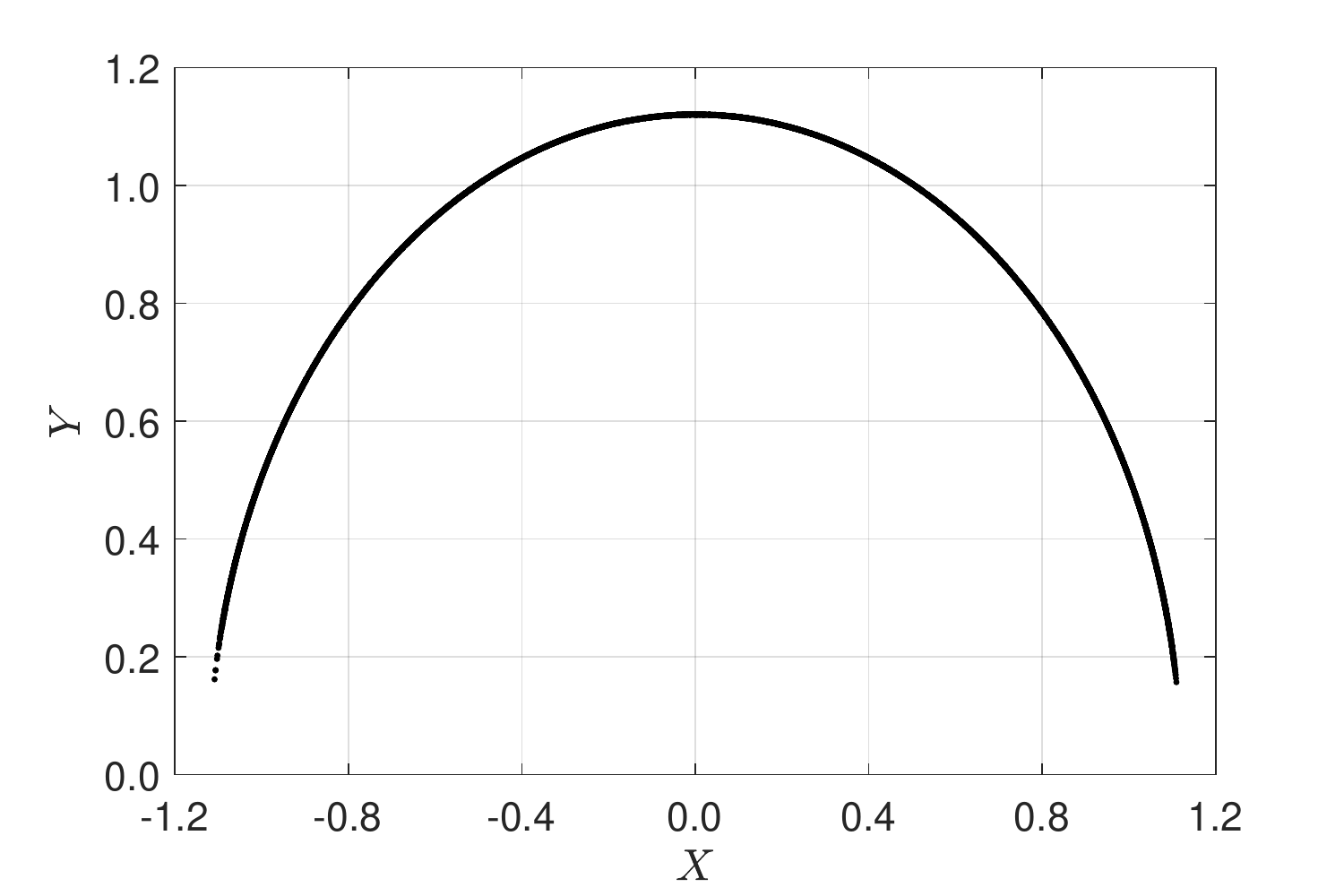}
		\caption{Nonlinear transformation of the slow variable $x$ onto a semicircle in the plane.}
		\label{fig:transform}
	\end{figure}

	To demonstrate the dimension reduction from a higher dimensional space, we can transform the slow variable $x$ by embedding it onto a curve in the plane (see Figure~\ref{fig:transform}).
	
\begin{figure}[htb]
		\centering
		\includegraphics[width=\linewidth]{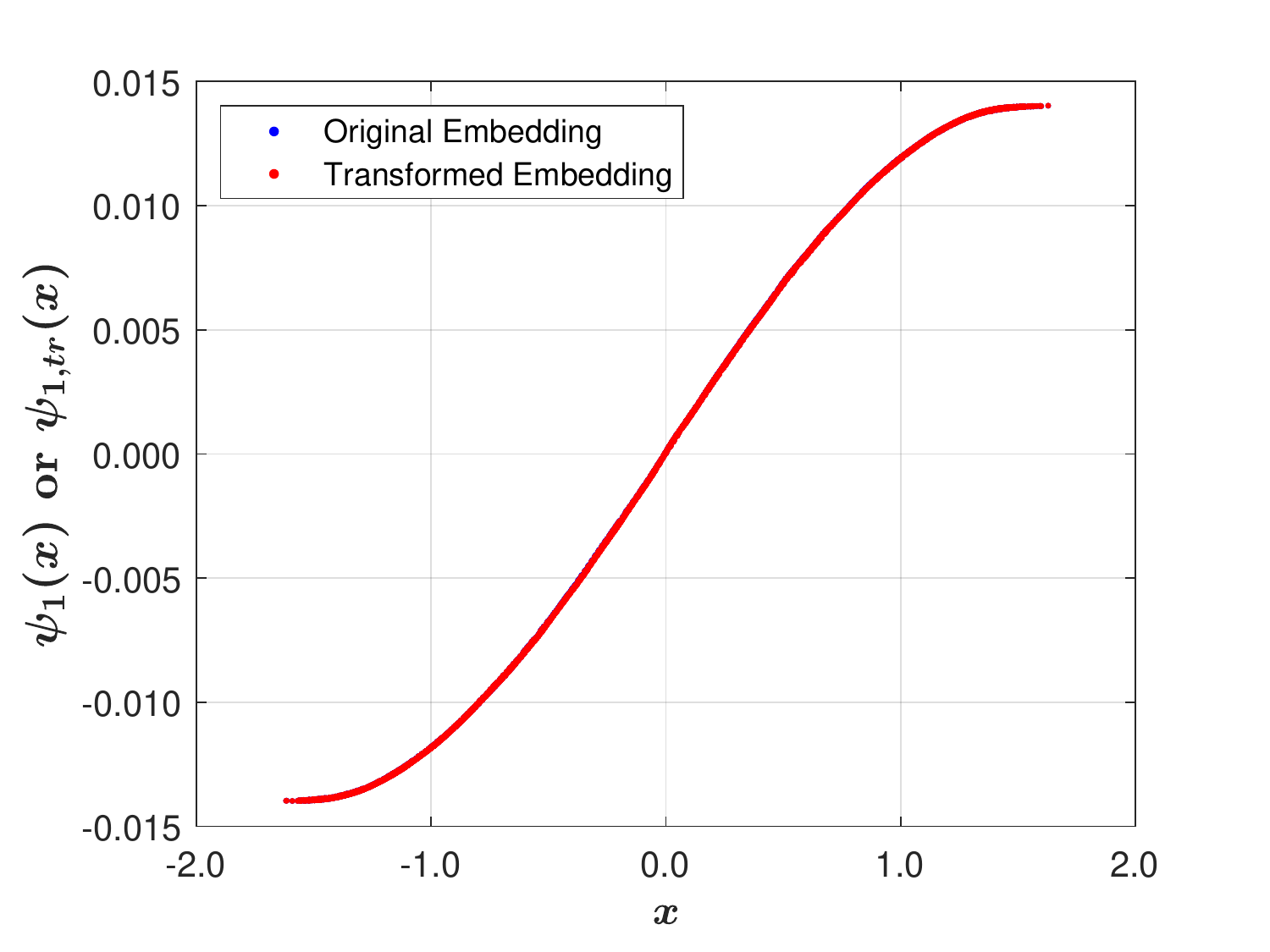}
		\caption{Embeddings of the original and transformed data set.}
		\label{fig:psi1Tr}
	\end{figure}

\begin{figure}[p]
		\centering
		\includegraphics[width=0.95\linewidth]{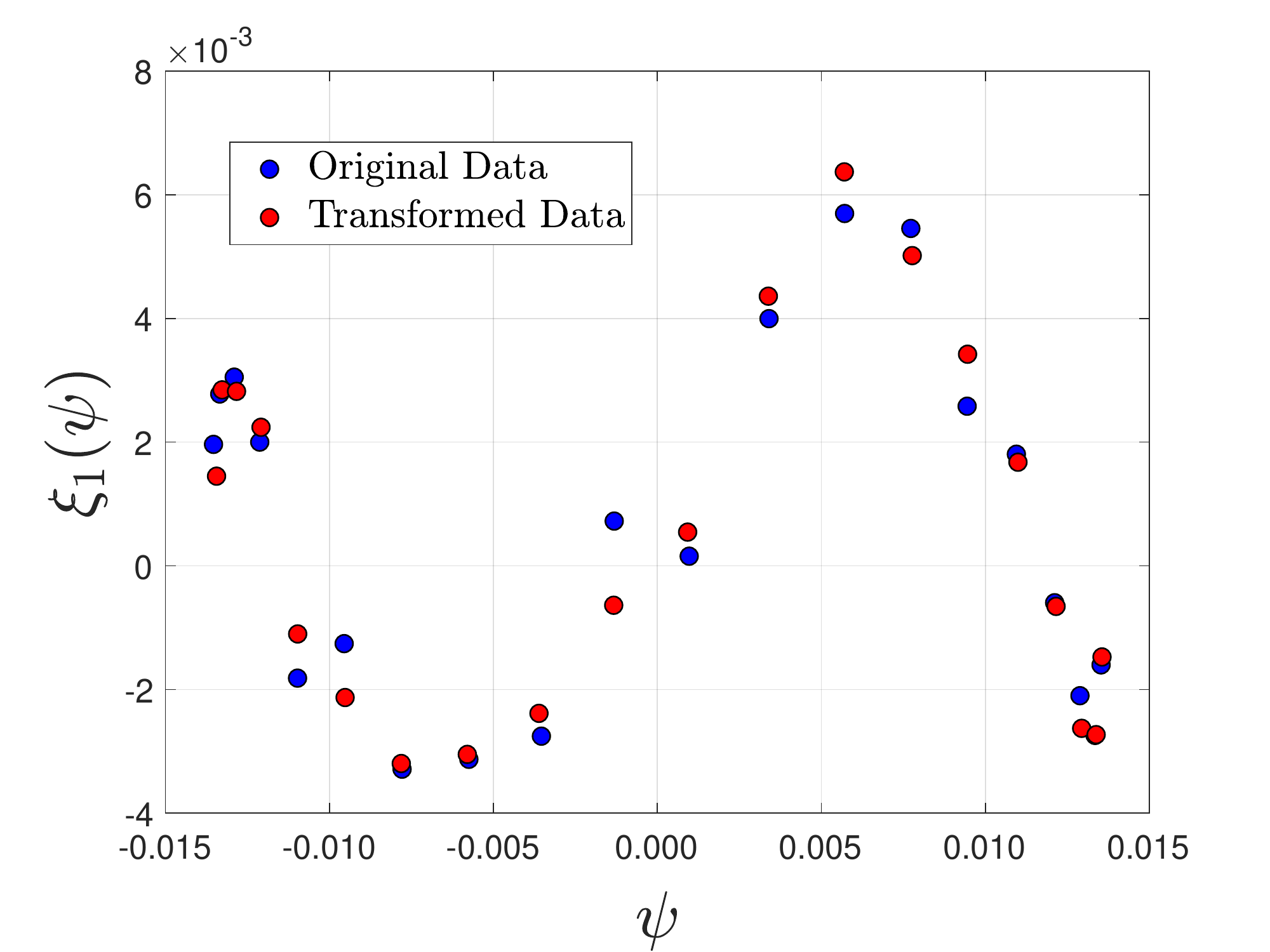}
		\caption{Estimation of the drift coefficient in Diffusion Map space.}
		\label{fig:ksi1Tr}
	\end{figure}

\begin{figure}[p]
		\centering
		\includegraphics[width=0.95\linewidth]{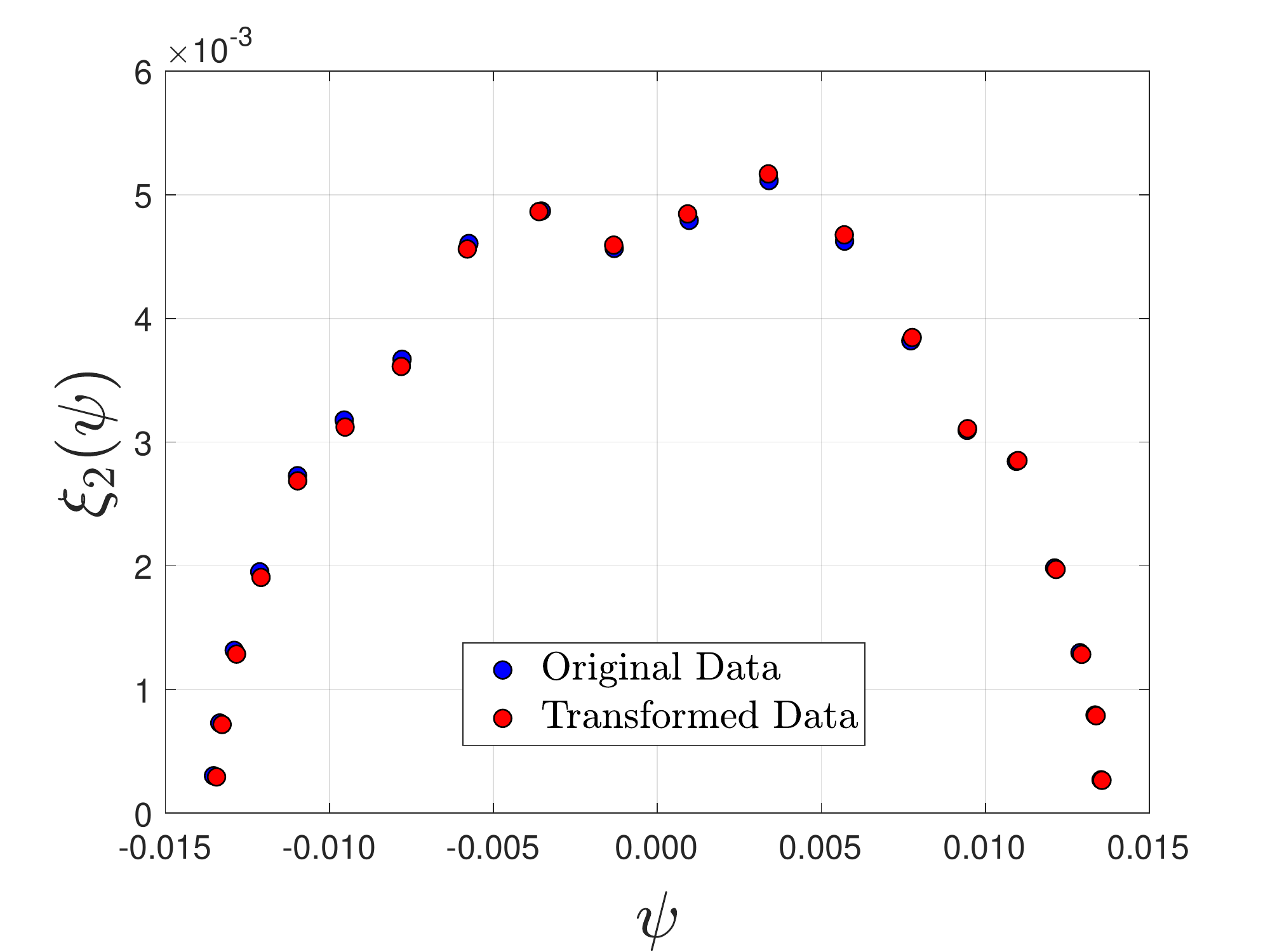}
		\caption{Estimation of the diffusion coefficient in Diffusion Map space.}
		\label{fig:ksi2Tr}
	\end{figure}

	The entire data set is now five-dimensional, and we can again apply DMaps with Mahalanobis distances as before. We again compute the pseudo-inverse using SVD as in~(\ref{eq:svd}), but now we use $d=n-1=4$ and discard the last singular value, which corresponds to the transverse direction on the semicircle. Figure~\ref{fig:psi1Tr} shows the embedding using the original data set as well as the transformed data set. The embeddings are almost identical. Using $\psi_{tr}(x)$, we can estimate again drift and diffusion coefficients as above. The results are shown in Figure~\ref{fig:ksi1Tr} and Figure~\ref{fig:ksi2Tr}.
	
	
	The same method can be applied in the case of multiplicative noise:
	\begin{equation*}
		\begin{aligned}
			\frac{dx}{dt} & = A(x-x^3)+\frac{\lambda}{\ep}(1+\nu x^2)y_2, & \qquad\frac{dy_1}{dt} & =\frac{10}{\ep^2}(y_2-y_1), \\
			\frac{dy_2}{dt} & =\frac{1}{\ep^2}(28y_1-y_2-y_1y_3), & \frac{dy_3}{dt} & =\frac{1}{\ep^2}(y_1y_2-\frac{8}{3}y_3).
		\end{aligned}
	\end{equation*}
	The approximate dynamics for the slow variable in this case are given by~\cite{krumscheid2013semiparametric}:
	\begin{equation*}
		dx = (Ax+Bx^3+Cx^5)\,dt+\sqrt{\sigma_a+\sigma_bx^2+\sigma_cx^4}\,dW.
	\end{equation*}
	The constants from~(\ref{eq:noiseSDE}) and the simulation parameters remain unchanged, with the exception of $\nu=1$ and $\Delta t=0.01$. Applying It\^{o}'s Lemma again, we obtain
	\begin{equation*}
		\begin{split}
			d\psi= & \left[(Ax+Bx^3+Cx^5)\frac{d\psi}{dx}+\frac{\sigma_a+\sigma_b x^2+\sigma_c x^4}{2}\frac{d^2\psi}{dx^2}\right]dt \\
			& +\sqrt{\sigma_a+\sigma_b x^2+\sigma_c x^4}\,\frac{d\psi}{dx}\,dW.
		\end{split}
	\end{equation*}
	
\begin{figure}[p]
		\centering
		\includegraphics[width=\linewidth]{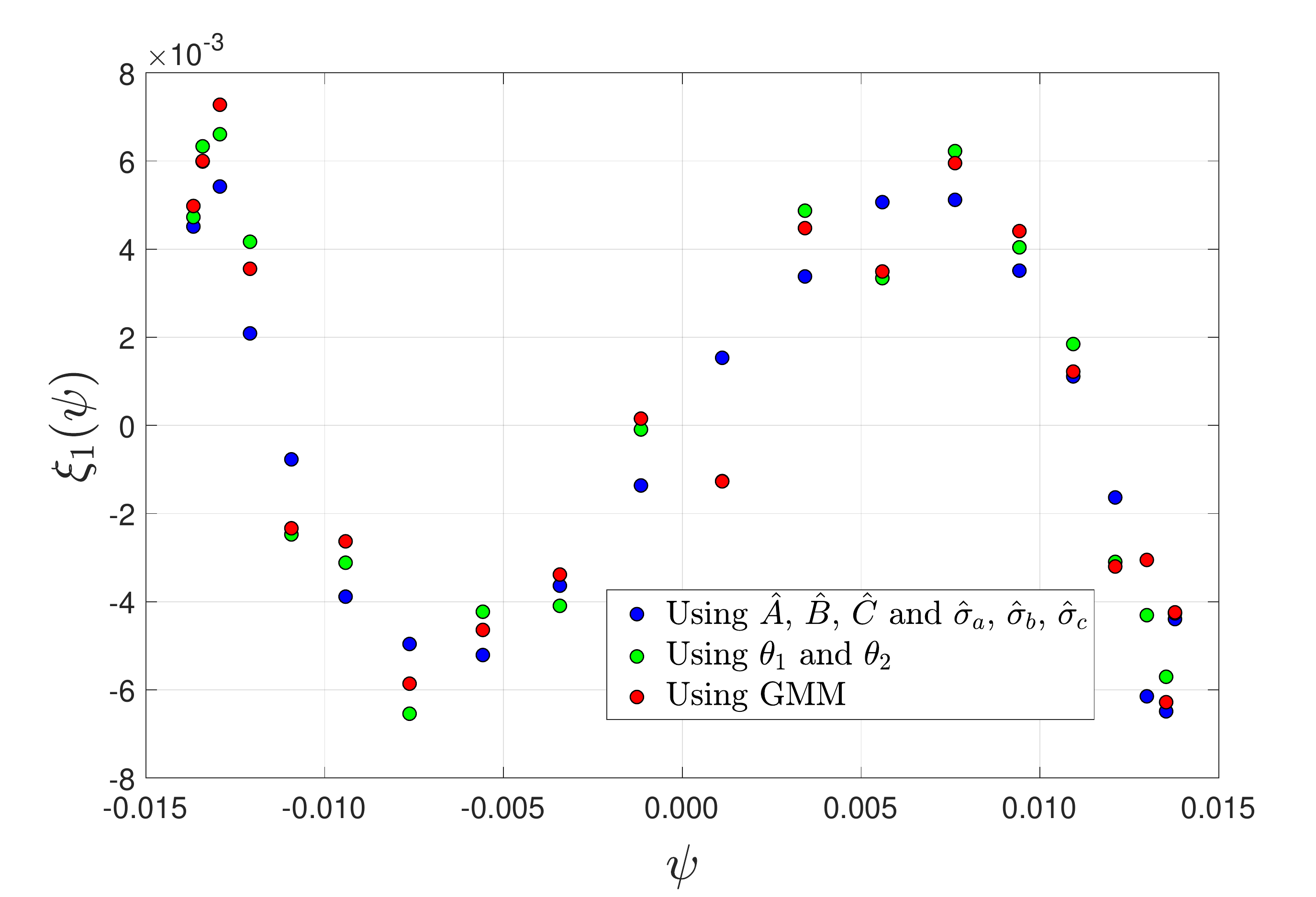}
		\caption{Estimation of the drift coefficient in diffusion map space.}
		\label{fig:ksi1_multi}
	\end{figure}

\begin{figure}[p]
		\centering
		\includegraphics[width=\linewidth]{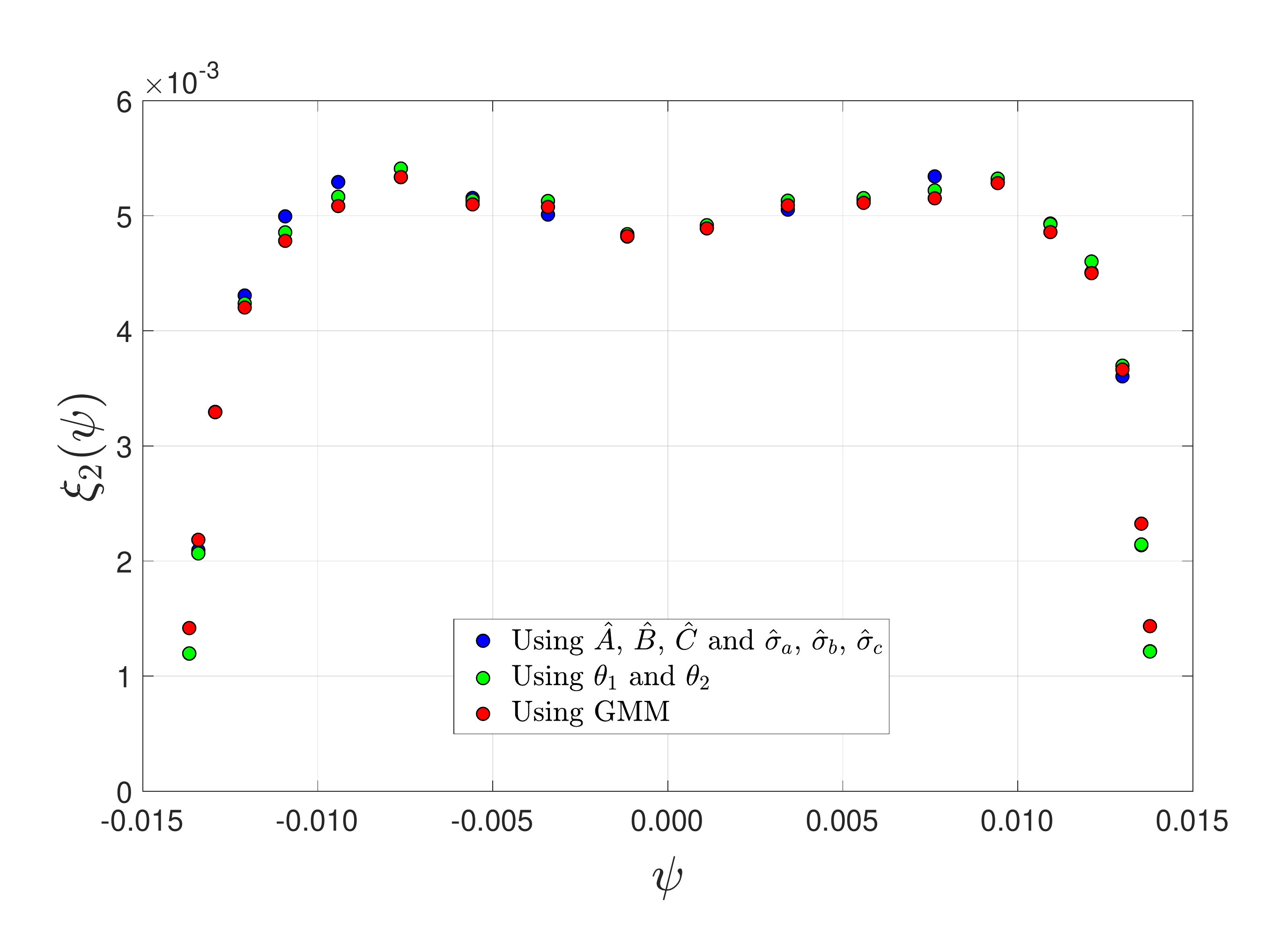}
		\caption{Estimation of the diffusion coefficient in diffusion map space.}
		\label{fig:ksi2_multi}
	\end{figure}

	If we estimate the coefficients of the polynomials in the drift and diffusion terms using the data in the original space,
	we can also estimate these coefficients in diffusion map space using an equation analogous to~(\ref{eq:1}):
	\begin{equation*}
		\begin{split}
			d\psi= & \left[(\hat{A}x+\hat{B}x^3+\hat{C}x^5)\frac{d\psi}{dx}+\frac{\hat{\sigma}_a+\hat{\sigma}_bx^2+\hat{\sigma}_cx^4}{2}\frac{d^2\psi}{dx^2}\right]dt \\
			& +\sqrt{\hat{\sigma}_a+\hat{\sigma}_b x^2+\hat{\sigma}_c x^4}\,\frac{d\psi}{dx}\,dW.
		\end{split}
	\end{equation*}
	The other two estimation methods remain the same, \textit{i.e.}, either using $\theta_1$ and $\theta_2$ calculated at each starting point or calculating $\xi_1$ and $\xi_2$ at each starting point. The results are shown in Figure~\ref{fig:ksi1_multi} and Figure~\ref{fig:ksi2_multi}.
	
	
	\section{Conclusions}
	We have confirmed that, at the limit of small time steps and large temperatures, trajectories produced by the Simulated Annealing / Random Walk Metropolis Hastings algorithm are analogous to those that result from the Langevin equation, a global, stochastic optimization algorithm. We use SA as our ``inner'' optimizer that produces ensembles of brief simulation bursts, which contain information about an effective gradient of the objective function. Using dimension reduction techniques such as PCA or, in the case of nonlinear manifolds, Diffusion Maps, we can obtain the parameterization of the underlying low-dimensional manifold.
	
	As our first example, we show that a two-dimensional, Bayesian model is effectively one-dimensional and DMaps can retrieve the important parameter (a sort of ``reaction coordinate'') for the optimization. Combining SA with DMaps achieves considerably faster approach to the maximum, compared with the simple SA alone.
	
	Starting from a two-dimensional SDE that corresponds to a Langevin equation for an objective function with two parameters, we derived the corresponding SDE in terms of the diffusion coordinates and approximated numerically the theoretical drift and diffusion coefficients. We then estimated the same drift and diffusion coefficients using parameter inference on the data set that came from the application of DMaps on the original trajectories. Additionally, we can transform the original data set through a nonlinear transformation by mapping it onto a portion of the surface of a sphere and apply DMaps on the transformed data set using Mahalanobis distances. In both cases the estimated parameters are closely comparable to the theoretical ones.
	
	For illustration purposes, we constructed a three-dimensional objective function that has a strongly attracting, two-dimensional manifold. This work constitutes a simple ``proof of concept'' acceleration demonstration for the classes of optimization problems we consider. It is also an illustration of the tools required to perform scientific computations (here, gradient descent) in a {\em latent} variable space, a space parameterized by on-the-fly processing of the data produced by the ``inner optimizer.'' Fast implementations of the techniques (like Geometric Harmonics) for translating back-and-forth between the original space, in which the problem was given, and the latent space, where the coarse optimization steps are taken, are crucial for the usefulness of the approach. The true benefits of this approach and its potential should be explored by applying it to truly high-dimensional problems where other methods slow down considerably. This is the subject of current work.
	
	\section*{Acknowledgments} This work was partially supported by the DARPA Lagrange program and an ARO MURI. The authors gratefully acknowledge many helpful discussions with Prof. Christian Kuehn of T. U. Munich, especially for  the case where the noise term arises through the coupling with chaotic dynamical systems.
	
	G.P.'s research is supported by the EPSRC, grants EP/P031587/1, EP/L024926/1, EP/L020564/1 and by JPMorgan Chase \& Co. Any views or opinions expressed herein are solely those of the authors listed, and may differ from the views and opinions expressed by JPMorgan Chase \& Co. or its affiliates. This material is not a product of the Research Department of J.P. Morgan Securities LLC. This material does not constitute a solicitation or offer in any jurisdiction.
	
	A.D. was supported by the Lloyd's Register Foundation Programme on Data Centric Engineering and by the Alan Turing Institute under the EPSRC grant [EP/N510129/1].
	
	\bibliographystyle{AIMS}
	\bibliography{refs}
\end{document}